\documentclass[10pt]{article}

\usepackage{amsmath}
\usepackage{amssymb}
\usepackage{indentfirst}
\usepackage{graphics} 
\usepackage{color}
\usepackage[utf8]{inputenc}

\makeatletter

\@addtoreset{equation}{section}
\makeatother
\def\<{\langle}
\def\>{\rangle}

\newtheorem{lem}{Lemma}[section]
\newtheorem{theo}{Theorem}[section]
\newtheorem{rem}{Remark}[section]

\makeatletter
   
   \@addtoreset{equation}{section}
\makeatother

\begin{document}
\title{\bf  Asymptotic profile of $L^{2}$-norm of solutions for wave equations with critical log-damping}

\author{Ruy Coimbra Char\~ao\thanks{Corresponding author: ruy.charao@ufsc.br}  \\{\small Department of Mathematics}, {\small Federal University of Santa Catarina} \\ {\small 88040-900, Florianopolis, SC,  Brazil,} 
\\
and\\Ryo Ikehata\thanks{ikehatar@hiroshima-u.ac.jp} \\ {\small Department of Mathematics}, \small{ Division of Educational Sciences}\\ {\small Graduate School of Humanities and Social Sciences} \\ {\small Hiroshima University} \\ {\small Higashi-Hiroshima 739-8524, Japan}}

\date{}
\maketitle
\begin{abstract}
We consider wave equations with a special type of log-fractional damping. We study the Cauchy problem for this model in ${\bf R}^{n}$, and we obtain an  asymptotic profile and optimal estimates of solutions as $t \to \infty$ in $L^{2}$-sense. A maximal discovery of this note is that under the effective damping, in case of $n=1$ $L^{2}$-norm of the solution blows up in infinite time, and in case of $n = 2$ $L^{2}$-norm of the solution never decays and never blows up in infinite time. The latter phenomenon seems to be a rare case in the community.   
\end{abstract}
\section{Introduction}
\footnote[0]{Keywords and Phrases: Wave equation; logarithmic damping; $L^{2}$-estimates; asymptotic profile; decay; growth.}
\footnote[0]{2010 Mathematics Subject Classification. Primary 35L05; Secondary 35B40, 35C20, 35S05.}

We consider the following dissipative wave equation:
\begin{align}
& u_{tt} -\Delta u + \mu Lu_t  = 0,\ \ \ (t,x)\in (0,\infty)\times {\bf R}^{n},\label{eqn}\\
& u(0,x)= u_0(x), \quad  u_{t}(0,x)= u_{1}(x),\ \ \ x\in{\bf R}^{n} ,\label{initial}
\end{align}
where $\mu > 0$ is a constant, and the linear operator  
\[L := \log(I + (-\Delta)^{\frac{1}{2}}): D(L) \subset L^{2}({\bf R}^{n}) \to L^{2}({\bf R}^{n}),\]
is defined by 
$$f \in D(L) := \{u \in L^{2}({\bf R}^{n})\,:\,\int_{{\bf R}^{n}}(1+\log^{2}(1+\vert\xi\vert))\vert\hat{u}(\xi)\vert^{2}d\xi < +\infty\},$$  
\[(Lf) (x) := {\cal F}_{\xi\to x}^{-1}\left(\log(1+\vert\xi\vert)\hat{f}(\xi)\right)(x).\]
\noindent
Here, one has just denoted the Fourier transform ${\cal F}_{x\to\xi}(f)(\xi)$ of $f(x)$ by 
\[{\cal F}_{x\to\xi}(f)(\xi) = \hat{f}(\xi) := \displaystyle{\int_{{\bf R}^{n}}}e^{-ix\cdot\xi}f(x)dx, \quad \xi \in {\bf R}^n,\]
as usual with $i := \sqrt{-1}$, and ${\cal F}_{\xi\to x}^{-1}$ expresses its inverse Fourier transform (for the logarithmic Laplacian itself including applications to PDEs, see \cite{B, B-1}, \cite{CW} and the references therein). 
\noindent
Concerning the existence of a unique solution to problem \eqref{eqn}-\eqref{initial}, by a similar argument to \cite[Proposition 2.1]{ITY} (see also \cite{CAI}) based on the Lumer-Phillips Theorem one can find that the problem (1.1)-(1.2) with initial data $(u_0, u_1) \in H^{1}({\bf R}^{n}) \times L^{2}({\bf R}^{n}) $  has a unique mild solution
\[u \in C([0,\infty);H^{1}({\bf R}^{n})) \cap C^{1}([0,\infty);L^{2}({\bf R}^{n})).\]

Let us first mention the motivation of this work by restricting to the following linear equations (mainly):
\begin{equation}\label{ikehata-1}
u_{tt}-\Delta u + \mu(-\Delta)^{\theta}u_{t} = 0,
\end{equation}
with $\mu > 0$ and $\theta \in [0,1]$. It is very well-known from many research articles produced so far in the case of $\theta \in[0,1/2)$ the solution to \eqref{ikehata-1} is diffusive (see e.g., \cite{DE}, \cite{DE-3}, \cite{DEP}, \cite{IT} ), while in the case when $\theta \in (1/2,1]$ the solution is oscillating (see e.g., \cite{CLI}, \cite{DEP}, \cite{Ike-14}, \cite{S}). However, it seems that there are few papers treating the special case of $\theta = 1/2$ in detail. For example, in \cite[Theorem 5]{DR} or \cite{NR} they derives $(L^{1}\cap L^{m})$-$L^{m}$ estimates with $m \in (1,2]$ such that
\[\Vert u(t,\cdot)\Vert_{m} \leq C(1+t)^{1-n(1-\frac{1}{m})}\Vert(u_{0},u_{1})\Vert_{L^{1}\cap L^{m}}.\]   
If we choose $m = 2$ one can get
\[\Vert u(t,\cdot)\Vert \leq C(1+t)^{1-\frac{n}{2}}\Vert(u_{0},u_{1})\Vert_{L^{1}\cap L^{2}}\]   
for all $n \geq 1$ and $\mu > 0$. This estimate is important to apply it to nonlinear problems, however, it seems unknown to check the optimality of those estimates. Indeed, the estimates for $n = 1$ (infinite time blowup)  and $n = 2$ (no-decay and no-blowup) may be trivial only from the viewpoint of the upper bound estimates. Optimality implies, in this research, both upper and lower bound estimates simultaneously.  This will be done by capturing the leading terms of solutions as time goes to infinity. In this connection, in the paper \cite[page 17]{DE-2} they speak of the case $\theta = 1/2$ as parabolic-like. As will be seen below, this fact depends on the value of $\mu > 0$ (this may be pointed out in \cite[Theorem 5]{DR}). Anyway, from the authors' point of views once we check previously published papers, one will see that the case $n = 1,2$ may be treated a little sloppy. It should be mentioned that our model (1.1) cannot be included as an example recenly studied in \cite{CI-0}, where several $L^{2}$-deay properties of the solution itself have been investigated to the wave equations with a more general damping. 
  
The purpose of this paper is to consider (1.1)-(1.2) by capturing asymptotic profiles (as $t \to \infty$), and is to study optimal rates of estimates for the solution itself in terms of $L^{2}$-norm. In particular, we shall deal with the low dimensional case $n = 1,2$ more politely. Furthermore, we are concentrating only on analyzing the $L^{2}$-norm itself of the solution, and this is because the higher order derivatives of solutions with respect to time and spatial variables do not reflect a singularity of them. A study on a singularity included in the solutions seems strongly attractive, so the analysis of the higher order derivatives is out of scope in our research.

In this work one sets 
\[P_{1} := \int_{{\bf R}^{n}}u_{1}(x)dx.\]

The result below is about the case of small $\mu$ satisfying $\mu \in (0,2)$. The solution has an oscillation property in asymptotic sense as $t \to \infty$.  
\begin{theo}\label{theorem-2}
Let $n \geq 1$ and $\mu \in (0,2)$. Let $u_{0} \in (H^{1}({\bf R}^{n})\cap L^{1}({\bf R}^{n}))$ and $u_{1} \in (L^{2}({\bf R}^{n})\cap L^{1,1}({\bf R}^{n}))$. Then, the solution $u(t,x)$ to problem \eqref{eqn}-\eqref{initial} satisfies
\[\Vert u(t,\cdot) - {\cal F}_{\xi\to x}^{-1}\left(\chi(t,\xi)\right)(\cdot)\Vert \leq  C\Big( \Vert u_{1}\Vert_{1,1} +  \Vert u_0\Vert_1\Big) \;t^{-\frac{n}{2}},\]
for $t \gg 1$, where $C > 0$ is a constant depending only on $n$ and $\mu \in (0,2)$, and
\[\chi(t,\xi) := (1+r)^{-\frac{t\mu}{2}}\frac{2}{r\sqrt{4-(\frac{\mu\log(1+r)}{r})^{2}}}\sin\left(\frac{r\sqrt{4-\mu^{2}}}{2}t\right)P_{1},\quad \xi \in {\bf R}^{n}, \,\,t\geq 0.\]
\end{theo}
\begin{theo}\label{theorem-2-1}
Under the same assumption as in Theorem {\rm \ref{theorem-2}} the solution $u(t,x)$ to problem {\rm (1.1)-(1.2)} satisfies
\[C_{1}\vert P_{1}\vert t^{1-\frac{n}{2}} \leq \Vert u(t,\cdot)\Vert 
\leq C_{2}\Big( \Vert u_{1}\Vert_{1,1} t^{-\frac{n}{2}} + \Vert u_{0}\Vert_{1} t^{-\frac{n}{2}}
+ |P_1| t^{1-\frac{n}{2}} + |P_1|e^{-\alpha t}
\Big),\]
for $t \gg 1$, where $C_{j} > 0$ {\rm (}$j = 1,2${\rm )} and $\alpha > 0$ are generous constants depending on $n$ and $\mu \in (0,2)$.
\end{theo}
The third result is dealing with the case of large $\mu > 2$. The oscillation property of the solution  in the zone of low frequency disappears  and have a diffusive aspect asymptotically. 

Let $\delta >0$ be a unique solution satisfying 
\begin{equation}\label{R-00}
\mu\log(1+\delta) - 2\delta = 0,
\end{equation}
$$\mu\log(1+r) - 2r > 0,\quad \forall r < \delta,$$
$$\mu\log(1+r) - 2r < 0,\quad \forall r > \delta.$$
A unique existence of such real number $\delta > 0$ can be guaranteed by a simple fundamental differential calculus. Furthermore, let
\[\lambda_{\pm} := - \frac{\mu \log(1+r)}{2} \pm \frac{\sqrt{\mu^{2}\log^{2}(1+r)-4r^{2}}}{2}, \quad \xi \in \{ r := |\xi| \leq \delta\}.\] 
Then, one can get the following statement.
\begin{theo}\label{theorem-3}
Let $n \geq 1$, $\mu > 2$ and let $\delta > 0$ be a constant defined by \eqref{R-00}. Suppose $[u_{0},u_{1}] \in (H^{1}({\bf R}^{n})\cap L^{1}({\bf R}^{n}))\times (L^{2}({\bf R}^{n})\cap L^{1,1}({\bf R}^{n}))$. Then, the solution $u(t,x)$ to problem \eqref{eqn}-\eqref{initial} satisfies
\[\Vert {\cal F}_{x \to \xi}(u(t,\cdot))(\xi) - \nu(t,\cdot)\Vert_{L^{2}(\vert\xi\vert\leq\delta)} \leq C(\Vert u_{0}\Vert + \Vert u_{0}\Vert_{1} + \Vert u_{1}\Vert + \Vert u_{1}\Vert_{1,1})t^{-\frac{n}{2}}, \quad t \gg 1,\]
and
\[\Vert {\cal F}_{x \to \xi}(u(t,\cdot))(\xi)\Vert_{L^{2}(\vert\xi\vert\geq\delta)} \leq C(\Vert u_{0}\Vert_{1} + \Vert u_{1}\Vert_{1})e^{-\alpha t}, \quad(t \gg 1),\]
where
\begin{equation}\label{ike-400}
\nu(t,\xi) := P_{1}\frac{e^{\lambda_{+}t} - e^{\lambda_{-}t} }{\sqrt{\mu^{2}\log^{2}(1+r)-4r^{2}}},\quad \forall \xi \in \{\vert\xi\vert < \delta\},
\end{equation}
where $C > 0$ is a constant depending only on $n$ and $\mu > 2$, and $\alpha > 0$ is a generous constant.
\end{theo}
\begin{rem}{\rm In some sense, $r = \delta$ is also a singular point in the frequency region $r := \vert\xi\vert \in [0,\infty)$. This is a peculiar phenomenon for the wave equation with log-damping satisfying $\mu > 2$. It seems that nobody has ever pointed out this new observation in the community.}
\end{rem}
As an application one can have the optimal estimate.
\begin{theo}\label{theorem-3-1}
Under the same assumption as in Theorem {\rm \ref{theorem-3}} the solution $u(t,x)$ to problem {\rm (1.1)-(1.2)} satisfies
\[C_{1}\vert P_{1}\vert t^{1-\frac{n}{2}} \leq \Vert u(t,\cdot)\Vert \leq C_{2}(\Vert u_{0}\Vert + \Vert u_{0}\Vert_{1} + \Vert u_{1}\Vert + \Vert u_{1}\Vert_{1,1}) t^{1-\frac{n}{2}}\]
for $t \gg 1$, where $C_{j} > 0$ {\rm (}$j = 1,2${\rm )} are generous constants depending on $n$ and $\mu \in (2,\infty)$.
\end{theo}
\begin{rem} {\rm In the case when $n = 1$ and $P_{1} \ne 0$ one notices from Theorems \ref{theorem-2-1} and \ref{theorem-3-1} that the $L^{2}$-norm of the solution $u(t,x)$ blows-up in infinite time with blowup rate $\sqrt{t}$, while in the two dimensional case the solution $u(t,x)$ never decays and never blows up. This phenomenon for $n = 2$ also seems quite rare case in the community (see \cite{CI} for a recent result). This type of study for all $\mu > 0$ seems to be unknown so far to the equation (1.1) with (at least) log-damping.}
\end{rem}
\begin{rem}
{\rm In particular, in the case when $\mu > 2$ from the profile point of view observed in Theorem \ref{theorem-3} the so called double diffusion structure is captured explicitly (see \eqref{ike-400} and Remark 4.1 below). This structure has been first discovered in D'Abbicco-Ebert \cite{DE} to the fractionally damped waves $u_{tt}-\Delta u + (-\Delta)^{\theta}u_{t} = 0$ with $0 < \theta < 1/2$, and developed in Piske-Chara\~o-Ikehata \cite{PCI} to the wave equation with a log-type damping
$u_{tt}-\Delta u + \log(I + (-\Delta)^{\theta})u_{t} = 0$ with $0 < \theta < 1/2$. A double diffusion structure for the case $\theta = 1/2$ and $\mu > 2$ does not seem to be clearly mentioned so far there.}
\end{rem}
Finally, let us study the critical case $\mu = 2$. One can explain a specialty of the critical case $\mu = 2$ by comparing already studied equation case (cf., \cite{D, DR}) such that
\begin{equation}\label{ike-300}
u_{tt}-\Delta u + \mu(-\Delta)^{1/2}u_{t} = 0.
\end{equation}
One observes that \eqref{ike-300} can be written  in the Fourier space as follows
\begin{equation}\label{ike-301}
\hat{u}_{tt} +\vert\xi\vert^{2}\hat{u} + \mu\vert\xi\vert\hat{u}_{t} = 0.
\end{equation}
In the case when $\mu = 2$ its corresponding characteristic roots for \eqref{ike-301} are $\lambda_{\pm} = -\frac{\vert\xi\vert}{2}$, that is, a diffusive aspect appears in the solution $\hat{u}(t,\xi)$, while, the characteristic roots for \eqref{eqn} with $\mu = 2$ are 
$$\lambda_{\pm} = -\frac{\log(1+\vert\xi\vert)}{2} \pm i\sqrt{\vert\xi\vert^{2}-\log^{2}(1+\vert\xi\vert)},$$
and since $r \geq \log(1+r)$ for all $r \geq 0$, one notices $\lambda_{\pm} \in {\bf C}$, which implies an oscillation property of the Fourier transformed solution $\hat{u}(t,\xi)$ of \eqref{eqn}.  There is a big difference between \eqref{ike-300} with $\mu = 2$ and \eqref{eqn} with $\mu = 2$. The following results express a peculiar property of the equation \eqref{eqn} with $\mu = 2$.
\begin{theo}\label{theorem-ike-1}
Let $n \geq 1$ and $\mu = 2$. Let $u_{0} \in (H^{1}({\bf R}^{n})\cap L^{1}({\bf R}^{n}))$ and  $u_{1} \in (L^{2}({\bf R}^{n})\cap L^{1,1}({\bf R}^{n}))$. Then, the solution $u(t,x)$ to problem \eqref{eqn}-\eqref{initial} satisfies
\[\Vert u(t,\cdot) - {\cal F}_{\xi\to x}^{-1}\left(\nu(t,\xi)\right)(\cdot)\Vert \leq  C\Big( \Vert u_{1}\Vert_{1,1} +  \Vert u_0\Vert_1\Big) \;t^{-\frac{n}{2}},\]
for $t \gg 1$, where $C > 0$ is a constant depending only on $n$, and
\[\nu(t,\xi) := (1+r)^{-t}\frac{1}{\sqrt{r^{2}-\log^{2}(1+r)}}\sin\left(t\sqrt{r^{2}-\log^{2}(1+r)}\right)P_{1},\quad \xi \in {\bf R}^{n}, \,\,t\geq 0.\]
\end{theo}
\begin{theo}\label{theorem-ike-2}
Under the same assumption as in Theorem {\rm \ref{theorem-ike-1}} the solution $u(t,x)$ to problem {\rm (1.1)-(1.2)} satisfies
\[C_{1}\vert P_{1}\vert t^{1-\frac{n}{2}} \leq \Vert u(t,\cdot)\Vert 
\leq C_{2}\Big( \Vert u_{1}\Vert_{1,1} t^{-\frac{n}{2}} + \Vert u_{0}\Vert_{1} t^{-\frac{n}{2}}
+ |P_1| t^{1-\frac{n}{2}} + |P_1|e^{-\alpha t}
\Big),\]
for $t \gg 1$, where $C_{j} > 0$ {\rm (}$j = 1,2${\rm )} and $\alpha > 0$ are generous constants depending on $n$.
\end{theo}
\begin{rem}{\rm As a result of a series of Theorems above, one can classify the asymptotic properties of the solution to problem (1.1)-(1.2) as follows: $0 < \mu \leq 2$ $\Rightarrow$ oscillating, while $2 < \mu$ $\Rightarrow$ diffusive. 
It seems interesting to note that although the structure of the solution for such cases is different, they have the same decay/nondecay rates for the $L^2$-norm of the solution itself. Furthermore, in the case of $n =1$, a strong singularity appears in the $L^{2}$-norm of the solution itself for all $\mu > 0$ which shows a growth property with the rate $\sqrt{t}$ as $t \to \infty$.}
\end{rem}

{\bf Notation.} {\small Throughout this paper, $\| \cdot\|_q$ stands for the usual $L^q({\bf R}^{n})$-norm. For simplicity of notation, in particular, we use $\| \cdot\|$ instead of $\| \cdot\|_2$. We also introduce the following weighted functional spaces.
\[L^{1,\gamma}({\bf R}^{n}) := \left\{f \in L^{1}({\bf R}^{n}) \; \bigm| \; \Vert f\Vert_{1,\gamma} := \int_{{\bf R}^{n}}(1+\vert x\vert^{\gamma})\vert f(x)\vert dx < +\infty\right\}.\]
For two positive functions $f(r)$ and $g(r)$ defined on $(0,\infty)$, one denotes $f(r) \sim g(r)$ as $r \to +0$ (resp. $r \to \infty$) if there exists constants $C_{j} > 0$ ($j = 1,2$) such that
\[C_{1}g(r) \leq f(r) \leq C_{2}g(r),\]
for small $r > 0$ (resp. $r \to \infty$). For two quantities $F(t)$ and $G(t)$ ($t \in [0,\infty)$), we denote $F(t) \preceq G(t)$ if there exists a constant $C > 0$ such that $F(t) \leq CG(t)$ for some $t \geq t_{0}$ with $t_{0} \geq 0$. Finally, we denote the surface area of the $n$-dimensional unit ball by $\omega_{n} := \displaystyle{\int_{\vert\omega\vert = 1}}d\omega$. }

The paper is organized as follows. In Section 2 we prove Theorems 1.1 and 1.5, and in Section 3 the proof of Theorem 1.3 is developed. Theorem 1.4 is proved in Section 4, which is rather technical. In Section 5 we finalize the proof of Theorem  1.2. A proof of Theorem \ref{theorem-ike-2} can be done in Section 6.


\section{Asymptotic profile: the case $0 < \mu \leq 2$}
In this section we deal with the case of $\mu \in (0,2]$, which corresponding to the non-effective damping case. Theorem \ref{theorem-2} with the case $\mu \in (0,2)$ and Theorem 1.5 for the case $\mu = 2$ can be treated simultaneously divided into two subsection, recpectively.\\ 

We first prepare the next useful lemmas which can be obtained from similar lemmas in \cite{CDI} and \cite{Log-damping} by a simple changing of variable.
\begin{lem}\label{lemma1}Let $p > -1$, $\delta >0$ and $\mu > 0$. Then, it holds that
\[\int_{0}^{\delta}(1+r)^{-\mu t}r^{p} dr \sim t^{-(p+1)} \quad (t \gg 1).\]
\end{lem}
\begin{lem}\label{lemma2}Let $p \in {\bf R}$,  $\delta >0$ and $\mu > 0$. Then, there exists a constant $\alpha > 0$ such that
\[\int_{\delta}^{\infty}(1+r)^{-\mu t}r^{p}dr \sim e^{-\alpha t} \quad (t \gg 1).\]
\end{lem}
In each one of  above two  lemmas the constants of equivalence may depend on $n, \mu, p$ and $\delta$.

We need to consider the equivalent problem to \eqref{eqn}--\eqref{initial} in the Fourier space  which is 
\begin{align}\label{F-eqn}
&w_{tt}(t,\xi) + |\xi|^2 w(t,\xi) + \mu\log(1+  |\xi|)w(t,\xi) = 0,\;\; (t,\xi)\in (0,\infty)\times {\bf R}^{n},\\
&w(0,\xi) = w_{0}(\xi), \quad
w_{t}(0,\xi) = w_1(\xi),  \;\; \xi\in  {\bf R}^{n}, \label{F-initial}
\end{align}
where
\[w(t,\xi) := \hat{u}(t,\xi), \quad w_{j}(\xi) := \hat{u}_{j}(\xi)\quad (j = 0,1).\]
\noindent
The characteristic roots associated to \eqref{F-eqn} are given by  
\begin{align}\label{roots}
\lambda_{\pm} = \dfrac{-\mu\log(1+|\xi|) \pm \sqrt{\mu^2\log^2(1+|\xi|) - 4| \xi|^{2}}}{2}, \quad \xi \in {\bf R}^n.
\end{align}

Now, we note that the function $f(r):=\mu \log(1+r)-2r$ satisfies $f(0)=0$ and $f(r) < 0 $ for all $r>0$ because $ f'(r)=\frac{\mu}{1+r} -2 < 0 $ for all $r  > 0$ due to the assumption  that $\mu \leq 2$. Therefore, the characteristic roots are  all complex-valued for all $\xi \in {\bf R}^n$, $\xi \ne 0$,  and are expressed by  
\begin{equation}\label{ike-9}
\lambda _{\pm} = - \frac{\mu \log(1+r)}{2} \pm \frac{r\sqrt{4-(\frac{\mu\log(1+r)}{r})^{2}}}{2}i,  \quad r=|\xi|, \;\; \xi \in {\bf R}^n . 
\end{equation}
Here, note that in the case when  $\mu \in (0,2]$ one has 
\[\mu\log(1+r)-2r \leq 0, \quad \forall r \geq 0.\]
Then, it is easy to check that the solution in the Fourier space can be explicitly expressed as
\[w(t,\xi)  = e^{-\frac{t\mu\log(1+r) }{2}}\cos\Big(\frac{r\sqrt{4-(\frac{\mu\log(1+r)}{r})^{2}}}{2}t \Big)w_{0}(\xi)\]
\[+e^{-\frac{t\mu \log(1+r) }{2}} \frac{\mu \log(1+r)}{r\sqrt{4-(\frac{\mu\log(1+r)}{r})^{2}}}\sin\Big(\frac{r\sqrt{4-(\frac{\mu\log(1+r)}{r})^{2}}}{2}t\Big)w_{0}(\xi)\]
\[+ e^{-\frac{t\mu}{2}\log(1+r)}\frac{2}{r\sqrt{4-(\frac{\mu\log(1+r)}{r})^{2}}}\sin\Big(\frac{r\sqrt{4-(\frac{\mu\log(1+r)}{r})^{2}}}{2}t\Big)w_{1}(\xi)\]
\[= (1+r)^{-\frac{t\mu  }{2}}\cos\Big(\frac{r\sqrt{4-(\frac{\mu\log(1+r)}{r})^{2}}}{2}t \Big)w_{0}(\xi)\]
\[+(1+r)^{-\frac{t\mu }{2}} \frac{\mu \log(1+r)}{r\sqrt{4-(\frac{\mu\log(1+r)}{r})^{2}}}\sin\Big(\frac{r\sqrt{4-(\frac{\mu\log(1+r)}{r})^{2}}}{2}t\Big)w_{0}(\xi)\]
\begin{equation}\label{ike-10}
+ (1+r)^{-\frac{t\mu}{2}}\frac{2}{r\sqrt{4-(\frac{\mu\log(1+r)}{r})^{2}}}\sin\Big(\frac{r\sqrt{4-(\frac{\mu\log(1+r)}{r})^{2}}}{2}t\Big)w_{1}(\xi).
\end{equation} 

By applying the decomposition of initial data $w_1(\xi)$  (see \cite{I-04}) such that
\begin{equation}\label{ike-1-1}
w_{1}(\xi) = \hat{u}_{1}(\xi) = P_{1} + A_{1}(\xi)-iB_{1}(\xi),
\end{equation}
where
\[P_{1} := \int_{{\bf R}^{n}}u_{1}(x)dx,\]
\[A_{1}(\xi) := \int_{{\bf R}^{n}}(\cos(x\xi)-1)u_{1}(x)dx, \quad B_{1}(\xi) := \int_{{\bf R}^{n}}\sin(x\xi)u_{1}(x)dx,\]
one can get the meaningful  identity   for $r=|\xi| >0$, with $\xi \in  {\bf R}^n$, as follows
\begin{align}\label{ike-11}
&w(t,\xi) - (1+r)^{-\frac{t\mu}{2}}\frac{2}{r\sqrt{4-(\frac{\mu\log(1+r)}{r})^{2}}}\sin\Big(\frac{r\sqrt{4-(\frac{\mu\log(1+r)}{r})^{2}}}{2}t\Big)P_{1} \nonumber\\
&= (1+r)^{-\frac{t\mu}{2}}\frac{2}{r\sqrt{4-(\frac{\mu\log(1+r)}{r})^{2}}}\sin\Big(\frac{r\sqrt{4-(\frac{\mu\log(1+r)}{r})^{2}}}{2}t\Big)(A_{1}(\xi)-iB_{1}(\xi)) \nonumber \\
& + (1+r)^{-\frac{t\mu }{2}} \frac{\mu \log(1+r)}{r\sqrt{4-(\frac{\mu\log(1+r)}{r})^{2}}}\sin\Big(\frac{r\sqrt{4-(\frac{\mu\log(1+r)}{r})^{2}}}{2}t\Big)w_0(\xi)\\
&+ (1+r)^{-\frac{t\mu  }{2}}\cos\Big(\frac{r\sqrt{4-(\frac{\mu\log(1+r)}{r})^{2}}}{2}t \Big)w_{0}(\xi). \nonumber
\end{align} 
\noindent
As a candidate of the leading term of the solution $w(t,\xi)$ as $t \to \infty$,  at this stage we can present the following function, defined for  $\xi \in {\bf R}^n\setminus\{0\}$ and $t \geq 0$,
\begin{eqnarray}\label{profile}
\label{ike-10-1}
\nu(t,\xi)  := (1+r)^{-\frac{t\mu}{2}}\frac{2}{r\sqrt{4-(\frac{\mu\log(1+r)}{r})^{2}}}\sin\Big(\frac{r\sqrt{4-(\frac{\mu\log(1+r)}{r})^{2}}}{2}t\Big)P_{1}.
\end{eqnarray}


\subsection{The case $ 0< \mu < 2$: Proof of Theorem 1.1}

 In this case it is possible to have the leading term in a slightly simpler form, applying the mean value theorem, and in fact, the profile $\nu(t,\xi)$ can be divided into the following form:
\begin{eqnarray}\label{profile1}
\nu(t,\xi) = (1+r)^{-\frac{\mu t}{2}}\frac{2}{r\sqrt{4-(\frac{\mu\log(1+r)}{r})^{2}}}\sin\Big(\frac{rt\sqrt{4-\mu^{2}}}{2}\Big)P_{1} + P_{1}R(t,\xi),
\end{eqnarray}
where
\begin{equation}\label{Remainder}
\hspace{-0.1cm} R(t,\xi) = \mu^{2}t (1+r)^{-\frac{\mu t}{2}}\frac{\cos(\eta(t,\xi))}{\sqrt{4-(\frac{\mu\log(1+r)}{r})^{2}}}\frac{1-\frac{\log^{2}(1+r)}{r^{2}}}{\sqrt{4-\mu^{2}} + \sqrt{4-(\frac{\mu\log(1+r)}{r})^{2}}},
\end{equation}
\[\eta(t,\xi) := \theta\frac{r\sqrt{4-(\frac{\mu\log(1+r)}{r})^{2}}}{2}t + (1-\theta)\frac{r\sqrt{4-\mu^{2}}}{2}t\]
with some $\theta \in (0,1)$.  
\noindent
Here, to simplify the notation one set 
\[\gamma := \frac{\sqrt{4-\mu^{2}}}{2} > 0.\]
Then, in case of $\mu \in (0,2)$ one can introduce a new simplified form of the leading term of the solution $w(t,\xi)$,  as $t \to \infty$, by 
\begin{equation}\label{Remainder2}
\chi(t,\xi) := (1+r)^{-\frac{\mu t}{2}}\frac{2}{r\sqrt{4-(\frac{\mu\log(1+r)}{r})^{2}}}\sin\Big(\gamma tr\Big)P_{1}.
\end{equation} 

Now, substituting the new expression \eqref{profile1} for $\nu(t,\xi)$   in the identity \eqref{ike-11}  and using the expression \eqref{Remainder2} for the simplified  leading term, one obtains the following  identity:
\begin{align}\label{ike-12}
&w(t,\xi) - \chi(t,\xi)
=w(t,\xi) - (1+r)^{-\frac{\mu t}{2}}\frac{2}{r\sqrt{4-(\frac{\mu\log(1+r)}{r})^{2}}}\sin\big(\gamma t r\big)P_{1} \nonumber\\
&= P_1R(t,\xi) +  (1+r)^{-\frac{t\mu}{2}}\frac{2}{r\sqrt{4-(\frac{\mu\log(1+r)}{r})^{2}}}\sin\Big(\frac{r\sqrt{4-(\frac{\mu\log(1+r)}{r})^{2}}}{2}t\Big)\big(A_{1}(\xi)-iB_{1}(\xi)\big) \nonumber \\
& + (1+r)^{-\frac{t\mu }{2}} \frac{\mu \log(1+r)}{r\sqrt{4-(\frac{\mu\log(1+r)}{r})^{2}}}\sin\Big(\frac{r\sqrt{4-(\frac{\mu\log(1+r)}{r})^{2}}}{2}t\Big)w_0(\xi)\\
&+ (1+r)^{-\frac{t\mu  }{2}}\cos\Big(\frac{r\sqrt{4-(\frac{\mu\log(1+r)}{r})^{2}}}{2}t \Big)w_{0}(\xi), \nonumber
\end{align} 
with the function $R(t,\xi)$ given by \eqref{Remainder}.

We have to estimate all terms in the right hand side of \eqref{ike-12}
including  $R(t,\xi)$.
 
In the following, let us make sure the validity of the fact that in the case when $\mu \in (0,2)$ the function given by \eqref{Remainder2} is really a leading term, by estimating the remainder terms given by functions $J_i(t,\xi), i=1,2,3, $ and $R(t,\xi)$ in terms of $L^{2}$-norm, where  $J_i(t,\xi)$ are three functions on the right hand side of \eqref{ike-12} defined as 
\begin{align*}
 J_1(t,\xi)& =: (1+r)^{-\frac{t\mu}{2}}\frac{2}{r\sqrt{4-(\frac{\mu\log(1+r)}{r})^{2}}}\sin\Big(\frac{r\sqrt{4-(\frac{\mu\log(1+r)}{r})^{2}}}{2}t\Big)\big(A_{1}(\xi)-iB_{1}(\xi)\big) ,\\
J_2(t,\xi) &
=(1+r)^{-\frac{t\mu }{2}} \frac{\mu \log(1+r)}{r\sqrt{4-(\frac{\mu \log(1+r)}{r})^{2}}}\sin\Big(\frac{r\sqrt{4-(\frac{\mu\log(1+r)}{r})^{2}}}{2}t\Big)w_0(\xi), \nonumber\\
J_3(t,\xi)&=(1+r)^{-\frac{t\mu  }{2}}\cos\Big(\frac{r\sqrt{4-(\frac{\mu\log(1+r)}{r})^{2}}}{2}t \Big)w_{0}(\xi). 
\end{align*}
\noindent
It is important to get upper bound estimates of the $L^2$-norm of each functions $J_i(t,\xi),\; i=1,2,3$ and $R(t,\xi)$.

For these  estimates one uses two inequalities such  as 
\[\left\vert \frac{\sin\theta}{\theta}\right\vert \leq K,\quad  \mbox{for all } \; \theta \in {\bf R},\]
with some $K \geq 1$, and
\begin{equation}\label{I-add-1}
\vert A_{1}(\xi)-iB_{1}(\xi)\vert \leq M\Vert u_{1}\Vert_{1,1}\vert\xi\vert,\quad \xi \in {\bf R}^{n},
\end{equation}
where $M > 0$ is a constant (see \cite{Ike-14}).
Then, by similar arguments as in \cite{Ike-14} (see also \cite{CAI}), by using \eqref{I-add-1} one can estimate the  $L^2$-norm of  $J_1(t,\xi)$ as follows:
\begin{align}\label{j1}
\int_{{\bf R}^{n}}\vert J_1(t,\xi)\vert^{2}d\xi
& \leq K^{2}M^{2}\Vert u_{1}\Vert_{1,1}^{2}t^{2}\int_{{\bf R}^{n}}(1+r)^{-\mu t}r^{2}d\xi \nonumber\\
&\leq K^{2}M^{2}\Vert u_{1}\Vert_{1,1}^{2}t^{2}\left(\int_{\vert\xi\vert\leq 1}(1+r)^{-\mu t}r^{2}d\xi + \int_{\vert\xi\vert\geq 1}(1+r)^{-\mu t}r^{2}d\xi\right) 
\\
&= K^{2}M^{2}\Vert u_{1}\Vert_{1,1}^{2}t^{2}\omega_{n}\left(\int_{0}^{1}(1+r)^{-\mu t}r^{n+1}dr + \int_{1}^{\infty}(1+r)^{-\mu t}r^{n+1}dr\right) \nonumber\\
&\leq \omega_{n}K^{2}M^{2}\Vert u_{1}\Vert_{1,1}^{2}t^{2}(t^{-n-2} + e^{-\alpha t}) \leq \omega_{n}K^{2}M^{2}\Vert u_{1}\Vert_{1,1}^{2}t^{-n} , \quad t \gg 1, \nonumber
\end{align}
where one has just used Lemmas \ref{lemma1} and \ref{lemma2} in the last line. 

While, one can get the estimate for $J_3(t,\xi)$ as follows,
\begin{align}\label{j3}
\int_{{\bf R}^{n}}\vert J_3(t,\xi)\vert^{2}d\xi & \leq \int_{{\bf R}^{n}}(1+r)^{-\mu t} \vert w_0(\xi)\vert^2d\xi \nonumber\\
& \leq ||u_0||_1^2 \Big( \int_{|\xi|\leq 1}(1+r)^{-\mu t}d\xi + \int_{|\xi| \geq 1}(1+r)^{-\mu t}d\xi\Big) \nonumber \\
&\leq C ||u_0||_1^2 \Big( t^{-n} + e^{-\alpha t} \Big)
\leq  C ||u_0||_1^2  t^{-n} , \quad t \gg 1,
\end{align}
with $C$ a positive constant depending on $\mu$ and $ n$. In the same way, using the fact $\mu\log(1+r) \leq  \mu r$ for all $r \geq 0$, and Lemmas \ref{lemma1} and \ref{lemma2} , we may get the following estimate for  $J_2(t,\xi)$,
\begin{align}
\label{j2}
\int_{{\bf R}^{n}}\vert J_2(t,\xi)\vert^{2}d\xi 
&\leq  
\int_{{\bf R}^{n}}(1+r)^{-\mu t}\mu^2 r^2 t^2 \vert w_0(\xi)\vert^2 d\xi \nonumber \\
&\leq  \mu^2 t^2\Vert u_{0}\Vert_{1}^{2}\;\big(t^{-n-2} + e^{-\alpha t}\big) \leq C \mu^{2}\Vert u_{0}\Vert_{1}^{2}\;t^{-n} , \quad t \gg 1 
\end{align}
with $C$ a positive constant.

Finally, one can estimate the part of the remainder terms given by the function $R(t,\xi)$ that is defined in \eqref{Remainder}. To do that we need the following trivial facts:
\begin{align}
&\lim_{r\rightarrow 0} \frac{r^2- \log^2(1+r)}{4r^2-\mu^2\log^2(1+r)}=1/4, \label{limit1}\\
&\lim_{r\rightarrow \infty} \frac{r^2- \log^2(1+r)}{4r^2-\mu^2\log^2(1+r)}=1/4.\label{limit2}
\end{align}

Moreover, one should note that for $0<\mu <2$ it holds that 
$$4-(\displaystyle{\frac{\mu\log(1+r)}{r}})^{2} > 0$$
for all $r > 0$. To simplify the notation as always one writes $r$ in place of $\vert \xi\vert$, thus one can estimate $R(t,\xi)$:
\begin{align}
&\int_{ \vert \xi \vert\leq 1} \vert R(t,\xi)\vert^2 d\xi \nonumber\\
 &  = \int_{\vert\xi\vert\leq 1}\Big|\mu^{2}t \;(1+r)^{-\frac{\mu t}{2}}\frac{\cos(\eta(t,\xi))}{\sqrt{4-(\frac{\mu\log(1+r)}{r})^{2}}}\frac{1-\frac{\log^{2}(1+r)}{r^{2}}}{\sqrt{4-\mu^{2}} + \sqrt{4-(\frac{\mu\log(1+r)}{r})^{2}}} \Big|^2d\xi\nonumber\\
  &  \leq \int_{ \vert \xi \vert\leq 1} \mu^{4}t^2 (1+r)^{-\mu t}\Big(\frac{1}{\sqrt{4-(\frac{\mu\log(1+r)}{r})^{2}}}\frac{1-\frac{\log^{2}(1+r)}{r^{2}}}{\sqrt{4-\mu^{2}} + \sqrt{4-(\frac{\mu\log(1+r)}{r})^{2}}} \Big)^2d\xi\nonumber\\
 &  \leq \int_{ \vert \xi \vert\leq 1}\frac{ \mu^{4}}{4-\mu^{2}}t^2 (1+r)^{-\mu t}\frac{\big(1-\frac{\log^{2}(1+r)}{r^{2}}\big)^2}{4-(\frac{\mu\log(1+r)}{r})^{2}} d\xi\nonumber\\
 &  = \int_{ \vert \xi \vert\leq 1}\frac{ \mu^{4}}{4-\mu^{2}}t^2 (1+r)^{-\mu t}\frac{\big(1-\frac{\log^{2}(1+r)}{r^{2}}\big)^2}{4-\frac{
 \mu^2\log^2(1+r)}{r^2}} d\xi\nonumber\\
 &= \int_{ \vert \xi \vert\leq 1}\frac{ \mu^{4}}{4-\mu^{2}}t^2 (1+r)^{-\mu t}\frac{\big(r^2-\log^{2}(1+r)\big)^2}{4r^2-
 \mu^2\log^2(1+r)} d\xi\nonumber\\
 &= \int_{ \vert \xi \vert\leq 1}\frac{ \mu^{4}}{4-\mu^{2}}t^2 (1+r)^{-\mu t}\frac{r^2-\log^{2}(1+r)}{4r^2-
 \mu^2\log^2(1+r)} \big(r^2-\log^{2}(1+r)\big) d\xi\nonumber\\
&\leq \int_{\vert \xi \vert\leq 1}\frac{ \mu^{4}}{4-\mu^{2}}t^2 (1+r)^{-\mu t}\frac{r^2-\log^{2}(1+r)}{4r^2-
 \mu^2\log^2(1+r)}r^2 d\xi. \nonumber
\end{align}

At this point, from \eqref{limit1} and \eqref{limit2} one can note that the function 
$$\frac{r^2-\log^{2}(1+r)}{4r^2 - \mu^2\log^2(1+r)}$$ 
is bounded on  $\xi \in {\bf R}^n$, $\xi \ne 0$. 
Hence, from the last estimate for $\int_{\vert\xi\vert\leq1}\vert R(t,\xi)\vert^2d\xi$ given by \eqref{Remainder} one obtains,
\begin{align}\label{Remainder3}
\int_{ \vert \xi \vert\leq 1} \vert P_1 R(t,\xi)\vert^2 d\xi 
& \leq P_1^2 C\frac{ \mu^{4}}{4-\mu^{2}} \int_{ \vert \xi \vert\leq 1}t^2 (1+r)^{-\mu t} r^2 d\xi  \nonumber \\
&= P_1^2\omega_n C\frac{ \mu^{4}}{4-\mu^{2}} \int_0^1t^2 (1+r)^{-\mu t} r^{n+1} d\xi \nonumber \\
& \leq  P_1^2C_{n,\mu} t^2t^{-n-2} \leq  P_1^2C_{n,\mu}t^{-n}, \quad t \gg 1,
 \end{align}
with $C_{n,\mu}>0$ a constant depending only on $\mu$ and $n$, where one has just used Lemma \ref{lemma1}.

However,  one can employ the same estimates as the low frequency zone to get the following estimates for $R(t,\xi)$ on the high frequency zone by relying on Lemma \ref{lemma2}, that is,
\begin{align}\label{Remainder4}
\int_{ \vert \xi \vert\geq 1} \vert P_1 R(t,\xi)\vert^2 d\xi 
& \leq P_1^2 C\frac{ \mu^{4}}{4-\mu^{2}} \int_{ \vert \xi \vert\geq 1}t^2 (1+r)^{-\mu t} r^2 d\xi  \nonumber \\
&= P_1^2\omega_n C\frac{ \mu^{4}}{4-\mu^{2}} \int_1^{\infty}t^2 (1+r)^{-\mu t} r^{n+1} d\xi \nonumber \\
& \leq  P_1^2C_{n,\mu} t^2 e^{-\alpha t} = P_1^2C_{n,\mu}  e^{-\frac{\alpha}{2} t}, \;\; \quad t \gg 1.
\end{align}

The  identity \eqref{ike-11} and estimates \eqref{j1}, \eqref{j3}, \eqref{j2}, \eqref{Remainder3}  and \eqref{Remainder4}  prove the next lemma.
\begin{lem} \label{ perfil-oscila}
Let $\mu \in (0,2)$, and let $\chi(t,\xi)$ be the function defined by \eqref{Remainder2}, and let $u_0 \in L^{1}({\bf R}^n)$ and  $u_{1} \in L^{1,1}({\bf R}^n)$. Then it is true that
\begin{align}\label{ profile-oscila}
\int_{ {\bf R}^n} |w(\xi,t) - \chi(\xi,t)|^2 d\xi  \leq C\Big(P_1^2+  \Vert u_{1}\Vert_{1,1}^{2} +  \Vert u_0\Vert^2_1\Big) \;t^{-n} + CP_{1}^{2}e^{-\alpha t}, \quad t \gg 1,
\end{align}
\end{lem}
with some constants $C=C(\mu,n) >0$ and $\alpha > 0$.

This implies the desired statement of Theorem \ref{theorem-2} which is polynomial decay.  

\subsection{The case $ \mu = 2$: Proof of Theorem 1.5}
In this case instead of using the simplified expression \eqref{Remainder2} for the asymptotic profile we use the leading term calculated by \eqref{ike-10-1} with $\mu = 2$, that is ,
\begin{eqnarray*}
\nu(t,\xi)  := (1+r)^{-\frac{t\mu}{2}}\frac{1}{r\sqrt{1-(\frac{\log(1+r)}{r})^{2}}}\sin\Big(r\sqrt{1-(\frac{\log(1+r)}{r})^{2}}t\Big)P_{1}
\end{eqnarray*}
defined for $t>0$ and $\xi  \in {\bf R}^n, \;\; \xi \ne 0.$ 

According to the estimates in Subsection 2.1, via \eqref{ike-11} the estimate for the $L^2$-norm of the difference $w(t,\xi) - \nu(t,\xi)$ depends only on the estimates for the functions $J_i(t,\xi), i=1,2,3$ already obtained 
in  \eqref{j1}, \eqref{j3} and \eqref{j2}. Thus, such fact proves Theorem 1.5.


\section{Asymptotic profile: the case $\mu > 2$}
In this section we deal with the case $\mu > 2$, which corresponding to the effective damping. 

For the case $\mu >2$ we note that the function $f(r):=\mu \log(1+r)-2r$  is such that $f(0)=0$ and $f'(r) = \displaystyle{\frac{\mu}{1+r}} -2 > 0 $ for  $0 < r  < \displaystyle{\frac{\mu -2}{2}}$ and $f$ has a positive global  maximum value in  $r= \displaystyle{\frac{\mu -2}{2}}$.  Then, one can conclude that there exists a unique number $\delta > \displaystyle{\frac{\mu -2}{2}}$ such that  the characteristic roots given by \eqref{roots}  satisfy 
\begin{align}\label{roots-real-com}
& \lambda _{\pm}  \mbox {are real for\;\;} 0 \leq r \leq \delta, \nonumber \\
&\lambda _{\pm}  \mbox {are  complex-valued  for}\;\;  r  > \delta.
\end{align}


\subsection{The case of real characteristic roots }
The characteristic  roots are real in the case $\mu>2$ and $r=|\xi| \leq \delta$.
 Here we  note that on the zone of low frequency $\{ \xi \in {\bf R}^n : |\xi| \leq \delta\}$ one has the following equivalence 
\begin{equation}\label{ike-608}
C_1 |\xi| \leq \mu \log(1+|\xi|) \leq  C_2 |\xi|, \quad |\xi| \leq \delta 
\end{equation}
 with $C_1$ and $C_2$ positive constants depending only on $\delta$.

In this case, the characteristic roots are real-valued given in \eqref{roots-real-com}, and using the above equivalence, they can be expressed as   
\begin{equation}\label{ike-19}
\lambda _{\pm} = - \frac{\mu \log(1+r)}{2} \pm \frac{\sqrt{\mu^{2}\log^{2}(1+r)-4r^{2}}}{2}, \quad \xi \in \{ r := |\xi| \leq \delta\}. 
\end{equation}
Then, it is easy to check that the solution of \eqref{F-eqn}--\eqref{F-initial} can be explicitly represented as
\begin{align}\label{ike-19-1}
w(t,\xi)& = \frac{e^{\lambda_{+}t}-e^{\lambda_{-}t}}{\lambda_{+}-\lambda_{-}}w_{1}(\xi) + \frac{\lambda_{+}e^{\lambda_{-}t}-\lambda_{-}e^{\lambda_{+}t}}{\lambda_{+}-\lambda_{-}}w_{0}(\xi) \nonumber\\
&=\frac{w_{1}(\xi)-\lambda_{-}w_{0}(\xi)}{\lambda_{+}-\lambda_{-}}e^{\lambda_{+}t} + \frac{w_{0}(\xi)\lambda_{+}-w_{1}(\xi)}{\lambda_{+}-\lambda_{-}}e^{\lambda_{-}t}, \quad |\xi| \leq \delta,
\end{align}
where $r := \vert\xi\vert$, $w(t,\xi) := \hat{u}(t,\xi)$ and $w_{j}(\xi) := \hat{u}_{j}(\xi)$ for $j = 0,1$ again. \\
\noindent
In the following, for simplicity of notation we set 
$$h(r) := \mu^{2}\log^{2}(1+r) - 4r^{2},$$
$L := \log(1+r)$, and $r:= \vert\xi\vert \leq \delta$. Then, basing on the explicit representation \eqref{ike-19-1} we use the Chill-Haraux idea \cite{CH} to have the relation:
\[w(t,\xi) = \frac{e^{-t\frac{(r^{2}+\lambda_{+}^{2})}{\mu L}}}{\lambda_{+}-\lambda_{-}}w_{1}(\xi) + e^{-\frac{tr^{2}}{\mu L}}w_{0}(\xi) + e^{-\frac{tr^{2}}{\mu L}}\frac{\lambda_{+}w_{0}(\xi)}{\lambda_{-}-\lambda_{+}}\]
\begin{equation}\label{ike-701}
+ e^{-\frac{tr^{2}}{\mu L}}\frac{\lambda_{-}(1- e^{-\frac{t\lambda_{+}^{2}}{\mu L}} )}{\lambda_{+}-\lambda_{-}}w_{0}(\xi) + e^{-\frac{tr^{2}}{\mu L}}\frac{w_{0}(\xi)\lambda_{+} - w_{1}(\xi)}{\lambda_{+}-\lambda_{-}} e^{-\frac{t\lambda_{-}^{2}}{\mu L}}.
\end{equation}
Let us apply the decomposition \eqref{ike-1-1} of the initial data $w_{1}(\xi)$ to \eqref{ike-701} to get the following equality on the zone of low frequency $|\xi| \leq \delta$ and $t > 0$. We use the fact
\[\lambda_{\pm}^{2} + r^{2} = -\mu L\lambda_{\pm} = -\mu\log(1+r)\lambda_{\pm}.\] 
Then, it follows that
\[w(t,\xi) = \frac{e^{-t\frac{(r^{2}+\lambda_{+}^{2})}{\mu L}}}{\lambda_{+}-\lambda_{-}}P_{1}  + \frac{e^{-t\frac{(r^{2}+\lambda_{+}^{2})}{\mu L}}}{\lambda_{+}-\lambda_{-}}(A_{1}(\xi)-iB_{1}(\xi)) + e^{-\frac{tr^{2}}{\mu L}}w_{0}(\xi) + e^{-\frac{tr^{2}}{\mu L}}\frac{\lambda_{+}w_{0}(\xi)}{\lambda_{-}-\lambda_{+}}\]
\[+ e^{-\frac{tr^{2}}{\mu L}}\frac{\lambda_{-}(1- e^{-\frac{t\lambda_{+}^{2}}{\mu L}} )}{\lambda_{+}-\lambda_{-}}w_{0}(\xi) + e^{-\frac{tr^{2}}{\mu L}}\frac{w_{0}(\xi)\lambda_{+} - w_{1}(\xi)}{\lambda_{+}-\lambda_{-}} e^{-\frac{t\lambda_{-}^{2}}{\mu L}}\]
\[=  \frac{P_{1}}{\lambda_{+}-\lambda_{-}}e^{t\lambda_{+}} -  \frac{P_{1}}{\lambda_{+}-\lambda_{-}}e^{t\lambda_{-}} + \frac{e^{-t\frac{(r^{2}+\lambda_{+}^{2})}{\mu L}}}{\lambda_{+}-\lambda_{-}}(A_{1}(\xi)-iB_{1}(\xi)) + e^{-\frac{tr^{2}}{\mu L}}w_{0}(\xi) + e^{-\frac{tr^{2}}{\mu L}}\frac{\lambda_{+}w_{0}(\xi)}{\lambda_{-}-\lambda_{+}}\]
\[+ e^{-\frac{tr^{2}}{\mu L}}\frac{\lambda_{-}(1- e^{-\frac{t\lambda_{+}^{2}}{\mu L}} )}{\lambda_{+}-\lambda_{-}}w_{0}(\xi) + \frac{\lambda_{+}}{\lambda_{+}-\lambda_{-}}e^{t\lambda{-}}w_{0}(\xi) - \frac{(A_{1}(\xi)-iB_{1}(\xi))}{\lambda_{+}-\lambda_{-}}e^{t\lambda_{-}}\]
\begin{equation}\label{ike-702}
=:\frac{P_{1}}{\lambda_{+}-\lambda_{-}}e^{t\lambda_{+}}  - \frac{P_{1}}{\lambda_{+}-\lambda_{-}}e^{t\lambda_{-}} + E_{1}(t,\xi) + E_{2}(t,\xi) + E_{3}(t,\xi) + E_{4}(t,\xi) + E_{5}(t,\xi) + E_{6}(t,\xi).
\end{equation}

Since $r^{2}+ \lambda_{+}^{2} = -\mu L\lambda_{+}$ and $\lambda_{+}-\lambda_{-} = \sqrt{h(r)}$, one has a meaningful equality such that
\begin{equation}\label{ike-704}
w(t,\xi) - P_{1}\frac{e^{t\lambda_{+}} -e^{t\lambda_{-}}  }{\sqrt{h(r)}} = \sum_{j=1}^{6}E_{j}(t,\xi).
\end{equation}
Now, let us introduce the leading term as $t \to \infty$ of the solution $w(t,\xi)$ in the low frequency region $r := \vert\xi\vert \leq \delta$ as follows:
\begin{equation}\label{ike-706}
\nu(t,\xi) := P_{1}\frac{e^{t\lambda_{+}}}{\sqrt{\mu^{2}\log^{2}(1+r)-4r^{2}}} - P_{1}\frac{e^{t\lambda_{-}}}{\sqrt{\mu^{2}\log^{2}(1+r)-4r^{2}}}, \quad \xi \in \{\xi \in {\bf R}^{n}\,:\,\vert\xi\vert \leq \delta\},\,\,t > 0.
\end{equation}
One prepares the following lemma to estimate such remainder terms \eqref{ike-704}.
\begin{lem}\label{ike-705} Let $\mu > 2$, and let $\delta > 0$ be a number defined in \eqref{roots-real-com}. Then, there exist real numbers $\delta_{1} \in (0,\delta]$, $c > 0$ and $d > 0$ such that \\
{\rm (1)}\,$c\log(1+r) \leq \lambda_{+} - \lambda_{-} = \sqrt{h(r)} \leq d\log(1+r)$ for $0 \leq r \leq \delta_{1}$,\\
{\rm (2)}\,$-cr \leq \lambda_{+} \leq -dr$ for $0 \leq r \leq \delta_{1}$,\\
{\rm (3)}\,$-c\log(+r) \leq \lambda_{-} \leq -d\log(1+r)$ for $0 \leq r \leq \delta_{1}$.
\end{lem}
{\it Proof of Lemma \ref{ike-705}.}\,\,We first prepare basic concept to prove results (1)-(3), and for later use.

Since
\[\lim_{r \to +0}\frac{\log(1+r)}{r} = 1,\]
there is a small number $\delta^{*} \leq \delta$ such that
\[\frac{1}{2}r \leq \log(1+r) \leq \frac{3}{2}r,\quad r \in [0,\delta^{*}].\]
One hand, the continuity of the function $g(r) := r^{-1}\log(1+r)$ on $[\delta^{*},\delta]$ guarantees the existence of two numbrs $M > m >0$ satisfying
\[m \leq g(r) \leq M,\quad r \in [\delta^{*},\delta].\]
Thus, one has the estimate
\begin{equation}\label{ike-902}
c_{1}r := \min\{\frac{1}{2},m\}r \leq \log(1+r) \leq \max\{\frac{3}{2},M\}r =: d_{1}r,\quad r \in [0,\delta].
\end{equation}

Under these preparations we shall prove the statements.

\underline{Proof of (1).} Let $\mu > 2$, an let $\delta_{1} \in (0,\delta]$ be a small number such that
\begin{equation}\label{ruy-100}
\mu^{2} - 4(1+\delta_{1})^{2} > 0.
\end{equation}
In fact, because of $\mu > 2$ one can choose as 
\[\delta_{1} \in (0,\min\{\frac{\mu}{2}-1,\delta\}].\]
\noindent
Then, 
\[\sqrt{h(r)} = \lambda_{+}-\lambda_{-} \geq \sqrt{\mu^{2}\log^2(1+r) - 4(1+\delta_{1})^{2}\log^{2}(1+r)}\]
\begin{equation}\label{ruy-101}
= \sqrt{\mu^{2}-4(1+\delta_{1})^{2}}\log(1+r),
\end{equation}
where one has just used the fact that
\begin{equation}\label{ruy-102}
r \leq (1+\delta_{1})\log(1+r), \quad \forall r \in [0,\delta_{1}].
\end{equation}
A check of \eqref{ruy-102} is an easy exercise of calculus. Also, it is trivial:
\begin{equation}\label{ike-901}
\sqrt{h(r)} \leq \mu\log(1+r),\quad r \in [0,\delta].
\end{equation} 
Therefore, from \eqref{ruy-101} and \eqref{ike-901} one can get the desired result (1) by choosing $c := \sqrt{\mu^{2}-4(1+\delta_{1})^{2}}$ and $d := \mu$:
\[d\log(1+r) \geq \sqrt{h(r)} \geq c\log(1+r),\quad r \in [0,\delta_{1}].\]

\underline{Proof of (2).} Indeed, from \eqref{ike-902} it follows that
\[-2\lambda_{+} = \frac{4r^{2}}{\mu\log(1+r) + \sqrt{h(r)}} \leq \frac{4r^{2}}{\mu}\frac{1}{\log(1+r)} \leq \frac{4}{c_{1}\mu}r, \quad r \in [0,\delta].\]
While, because of (1), \eqref{ike-901} and \eqref{ike-902} one has
\[-2\lambda_{+} = \frac{4r^{2}}{\mu\log(1+r) + \sqrt{h(r)}} \geq \frac{4r^{2}}{2\mu\log(1+r)} \geq \frac{2r^{2}}{\mu d_{1}r} =: d_{2}r , \quad r \in [0,\delta_{1}].\]
These prove the statement of (2) in $[0,\delta_{1}]$.

\underline{Proof of (3).} Indeed, from \eqref{ike-901} one has 
\[2\lambda_{-} = -\mu\log(1+r) - \sqrt{h(r)} \geq -2\mu \log(1+r).\]
While, $2\lambda_{-} \leq -\mu\log(1+r)$ is trivial. These mean the statement of (3) in $[0,\delta]$. 

It should be remarked, in fact, one can choose the same coefficients $c >0$ and $d > 0$ through all results (1)-(3) by a standard argument.
\hfill
$\Box$

Based on Lemma \ref{ike-705} let us check that each terms $E_{j}(t,\xi)$ ($j = 1,2,3,4,5$) are remainders in the decomposition \eqref{ike-704}.


\underline{(1)\,Low frequency estimate for $w(t,\xi)$ on $[0,\delta_{1}]$ when $\mu > 2$.}

Let $\delta_{1} > 0$ be the number defined in Lemma \ref{ike-705}. In this part we estimate each remainder terms of \eqref{ike-704} on $r \in [0,\delta_{1}]$ becaue we use Lmma \ref{ike-705} to check them. 

At first one can estimate $E_{1}(t,\xi)$ as follows by using (1) and (2) of Lemma \ref{ike-705}, the decomposition \eqref{ike-1-1} of the initial data and \eqref{I-add-1}:
\[\int_{\vert\xi\vert\leq\delta_{1}}\vert E_{1}(t,\xi)\vert^{2}\xi \leq \int_{\vert\xi\vert\leq\delta_{1}}\frac{e^{2t\lambda_{+}}}{\vert \lambda_{+}-\lambda_{-}\vert^{2}}\vert A_{1}(\xi) - iB_{1}(\xi)\vert^{2}d\xi\]
\[\preceq M^{2}\Vert u_{1}\Vert_{1,1}^{2}\int_{\vert\xi\vert\leq\delta_{1}}\frac{r^{2}e^{-ctr}}{\log^{2}(1+r)}d\xi \preceq \Vert u_{1}\Vert_{1,1}^{2}\int_{0}^{\delta_{1}}e^{-ctr}r^{n-1}dr\]
\begin{equation}\label{ike-800}
\leq C\Vert u_{1}\Vert_{1,1}^{2}t^{-n}, \quad t > 0
\end{equation}
with some equivalence constant $c > 0$.

Next, let us check $E_{2}(t,\xi)$, however, this is easy to have the estimate:
\begin{equation}\label{ike-801}
\int_{\vert\xi\vert\leq\delta_{1}}\vert E_{2}(t,\xi)\vert^{2}\xi \leq C\Vert u_{0}\Vert_{1}^{2}t^{-n}, \quad t > 0
\end{equation}
because of \eqref{ike-902}.

For $E_{3}(t,\xi)$, by using (1) and (2) of Lemma \ref{ike-705} and \eqref{ike-902} one has the following:
\[\int_{\vert\xi\vert\leq\delta_{1}}\vert E_{3}(t,\xi)\vert^{2}\xi \leq \int_{\vert\xi\vert\leq\delta_{1}} e^{-\frac{2tr^{2}}{\mu\log(1+r)}}\frac{\lambda_{+}^{2}}{\vert \lambda_{-}-\lambda_{+} \vert^{2}}\vert w_{0}(\xi)\vert^{2}d\xi\]
\[\preceq \Vert u_{0}\Vert_{1}^{2}\int_{\vert\xi\vert\leq\delta_{1}} e^{-\frac{2tr^{2}}{\mu\log(1+r)}}\frac{r^{2}}{\log^{2}(1+r)}d\xi \preceq \Vert u_{0}\Vert_{1}^{2}\int_{0}^{\delta_{1}}e^{-ctr}r^{n-1}dr\]
\begin{equation}\label{ike-802}
\leq C\Vert u_{0}\Vert_{1}^{2}t^{-n}, \quad t > 0
\end{equation}
with some equivalence constant $c > 0$.

$E_{4}(t,\xi)$ can be estimated by using (1) and (3) of Lemma \ref{ike-705} and \eqref{ike-902}:
\[\int_{\vert\xi\vert\leq\delta_{1}}\vert E_{4}(t,\xi)\vert^{2}\xi \leq \int_{\vert\xi\vert\leq\delta_{1}} e^{-\frac{2tr^{2}}{\mu\log(1+r)}}\frac{\lambda_{-}^{2}}{\vert \lambda_{+}-\lambda_{-} \vert^{2}}\vert 1- e^{-\frac{t\lambda_{+}^{2}}{\mu\log(1+r)}}\vert^{2}\vert w_{0}(\xi)\vert^{2}d\xi\]
\[\preceq \Vert u_{0}\Vert_{1}^{2}\int_{\vert\xi\vert\leq\delta_{1}} e^{-ctr}d\xi \]
\begin{equation}\label{ike-803}
\leq C\Vert u_{0}\Vert_{1}^{2}t^{-n}, \quad t > 0
\end{equation}
with some equivalence constant $c > 0$.

We treat $E_{5}(t,\xi)$ by using  can be estimated by Lemma \ref{ike-705} and Lemma 2.1 :
\[\int_{\vert\xi\vert\leq\delta_{1}}\vert E_{5}(t,\xi)\vert^{2}d\xi \leq \int_{\vert\xi\vert\leq\delta_{1}} e^{2t\lambda_{-}}\frac{\vert \lambda_{+}\vert^2}{\vert \lambda_{+}-\lambda_{-}\vert^{2}}\vert w_{0}(\xi)\vert^{2} d\xi\]
\[\preceq \Vert u_{0}\Vert_{1}^{2}\int_{\vert\xi\vert\leq\delta_{1}} e^{-ct\log(1+r)}d\xi \]
\begin{equation}\label{ike-803}
\leq C\Vert u_{0}\Vert_{1}^{2}t^{-n}, \quad t > 0
\end{equation}
with some equivalence constant $c > 0$.

We treat $E_{6}(t,\xi)$ can be checked by \eqref{ike-1-1}, \eqref{I-add-1}, Lemma \ref{ike-705} and Lemma 2.1 :
\[\int_{\vert\xi\vert\leq\delta_{1}}\vert E_{6}(t,\xi)\vert^{2}d\xi \leq \int_{\vert\xi\vert\leq\delta_{1}} e^{2t\lambda_{-}}\frac{\vert A_{1}(\xi)-iB_{1}(\xi)\vert^2}{\vert \lambda_{+}-\lambda_{-}\vert^{2}}d\xi\]
\[\preceq M^{2}\Vert u_{1}\Vert_{1,1}^{2}\int_{\vert\xi\vert\leq\delta_{1}} e^{2t\lambda_{-}}\frac{r^{2}}{\vert \lambda_{+}-\lambda_{-} \vert^{2}}d\xi \]
\[\preceq M^{2}\Vert u_{1}\Vert_{1,1}^{2}\int_{\vert\xi\vert\leq\delta_{1}} e^{-ct\log(1+r)}\frac{r^{2}}{\log^{2}(1+r)}d\xi \]
\begin{equation}\label{ike-803}
\leq C\Vert u_{1}\Vert_{1,1}^{2}t^{-n}, \quad t > 0
\end{equation}
with some equivalence constant $c > 0$.

Summarizing computations above one can arrive at the following low frequency estimate.
\begin{lem}\label{lem4.1}
Let $n \geq 1$ and $\mu > 2$. Then, it holds that
\[\int_{|\xi| \leq \delta_{1}} |w(t,\xi) - \nu(t,\xi)|^2d\xi \leq C_{n,\mu} \big(\|u_0\|_{1}^2 + \|u_1\|_{1,1}^2 \big)\; t^{-n}, \quad t \gg 1,\]
where $C=C_{n,\mu}$ is a positive constant and $\beta >0$ is a constant.
\end{lem}

Next we deal with th remainder terms of \eqref{ike-704} in the middle frequency region $[\delta_{1},\delta]$\\

\underline{(2)\,Midle frequency estimate for $w(t,\xi)$ on $[\delta_{1},\delta]$ when $\mu > 2$.}

To prove the middle frequency estimate for the solution $w(t,\xi)$ we rely on the energy estimate in the Fourier space due to \cite[Proposition 2.1]{CI-0}:
\[\vert\xi\vert^{2}\vert w(t,\xi)\vert^{2} \leq Ce^{-c\rho(\vert\xi\vert)t}(\vert\xi\vert^{2}\vert w_{0}(\xi)\vert^{2} + \vert w_{1}(\xi)\vert^{2}), \quad \xi \in {\bf R}^{n}\setminus\{0\},\]
where
\[\rho(r) := \frac{\mu r^{2}\log(1+r)}{r^{2}+\mu^{2}\log^{2}(1+r)}.\]
Note that Proposition 2.1 of \cite{CI-0} itself can be applied to our case with the weight function $\rho(r)$ just defined. This implies
\[\vert w(t,\xi)\vert^{2} \leq Ce^{-c\rho(\vert\xi\vert)t}(\vert w_{0}(\xi)\vert^{2} + \vert\xi\vert^{-2}\vert w_{1}(\xi)\vert^{2}), \quad\xi \in {\bf R}^{n}\setminus\{0\}.\]
Thus, for $r \in [\delta_{1},\delta]$ one has
\[\vert w(t,\xi)\vert^{2} \leq Ce^{-c\rho(\vert\xi\vert)t}(\vert w_{0}(\xi)\vert^{2} + \delta_{1}^{-2}\vert w_{1}(\xi)\vert^{2}) \leq C(1+\delta_{1}^{-2})e^{-c\rho(\vert\xi\vert)t}(\vert w_{0}(\xi)\vert^{2} + \vert w_{1}(\xi)\vert^{2}).\] 
This implies
\[\int_{\delta_{1}\leq\vert\xi\vert\leq \delta}\vert w(t,\xi)\vert^{2}d\xi \leq C(1+\delta_{1}^{-2})\int_{\delta_{1}\leq\vert\xi\vert\leq \delta}e^{-c\rho(\vert\xi\vert)t}(\vert w_{0}(\xi)\vert^{2} + \vert w_{1}(\xi)\vert^{2})d\xi.\]
In the case when  $r \in [\delta_{1},\delta]$ one has
\[e^{-c\rho(\vert\xi\vert)t} \leq e^{-Act},\]
where
\[A:= \frac{\mu\delta_{1}^{2}\log(1+\delta_{1})}{\delta^{2}+\mu^{2}\log^{2}(1+\delta)}.\]
Therefore, it follows that
\begin{equation}\label{ike-955}
\int_{\delta_{1}\leq\vert\xi\vert\leq \delta}\vert w(t,\xi)\vert^{2}d\xi \leq C(1+\delta_{1}^{-2})e^{-Act}(\Vert u_{0}\Vert^{2} + \Vert u_{1}\Vert^{2}),
\end{equation}
which implies the exponential decay result of the solution $w(t,\xi)$ in $[\delta_{1},\delta]$.\\

Next, let us derive exponential decay for the leading term $\nu(t,\xi)$ itself in $[\delta_{1},\delta]$.\\
For this aim set
\[A^{*} := \frac{\delta_{1}^{2}}{\mu\log(1+\delta) + \sqrt{\mu^{2}\log^{2}(1+\delta)-4\delta_{1}^{2}}}.\]
Note that since $\delta_{1} < \delta$, one has $\mu\log(1+\delta) > 2\delta_{1}$ which implies the well-definedness of $A^{*}$. Furthermore, it should be noted that
since
\[\lim_{r \to +0}\frac{1-e^{-r}}{r} = 1,\]
as in the argument of Lemma \ref{ike-705} one can assume that there exists constants $\gamma_{0} > 0$ and $\delta_{2} > 0$ such that
\begin{equation}\label{CR-3}
\gamma_{0} \leq \frac{1-e^{-r}}{r},\quad \forall r \in (0,\delta_{2}].
\end{equation}
In particular, since
\[\lim_{r \to \infty}\frac{1-e^{-r}}{r} = 0,\]
and the continuity of the function $r \mapsto \frac{1-e^{-r}}{r}$ on $(0,\infty)$, there exists a number $\kappa_{0} > 0$ such that
\begin{equation}\label{CR-7}
0 < \frac{1-e^{-r}}{r} \leq \kappa_{0},\quad \forall r > 0.
\end{equation}
Then, it follows from the general above result \eqref{CR-7}, one has
\[\frac{\vert 1-e^{-(\lambda_{+}-\lambda_{-})t} \vert^{2}}{\vert \lambda_{+}-\lambda_{-}\vert^{2}t^{2}} \leq \kappa_{0}^{2},\quad t > 0,\quad r \in [\delta_{1},\delta].\]
Thus, one can estimate as follows:
\[K_{0}(t) := \int_{\delta_{1}\leq\vert\xi\vert < \delta}\vert\nu(t,\xi)\vert^{2}d\xi = P_{1}^{2}t^{2}\int_{\delta_{1}\leq\vert\xi\vert < \delta}e^{2\lambda_{+}t}\frac{\vert 1-e^{-(\lambda_{+}-\lambda_{-})t} \vert^{2}}{\vert \lambda_{+}-\lambda_{-}\vert^{2}t^{2}} d\xi\]
\[\leq \kappa_{0}^{2}P_{1}^{2}t^{2}\int_{\delta_{1}\leq\vert\xi\vert < \delta}e^{2\lambda_{+}t}d\xi = \kappa_{0}^{2}P_{1}^{2}t^{2}\omega_{n}\int_{\delta_{1}}^{\delta}e^{\frac{-4t r^{2}}{\mu\log(1+r) + \sqrt{\mu^{2}\log^{2}(1+r)-4r^{2}}}}r^{n-1}dr\]
\[\leq \kappa_{0}^{2}P_{1}^{2}t^{2}\omega_{n}\int_{\delta_{1}}^{\delta}e^{-4A^{*}t}r^{n-1}dr \leq \kappa_{0}^{2}P_{1}^{2}\omega_{n}t^{2}e^{-4A^{*}t}\delta^{n-1}(\delta-\delta_{1}),\]
which implies
\begin{equation}\label{ike-951}
K_{0}(t) \leq CP_{1}^{2} e^{-\alpha t}\quad (t > 1).
\end{equation}
Now, let us prove the first estimate of Theorem 1.3. It follows from Lemma \ref{lem4.1}, \eqref{ike-955} and \eqref{ike-951} that
\[\Vert w(t,\cdot)-\nu(t,\cdot)\Vert_{L^{2}(\vert\xi\vert\leq\delta)} = \Vert w(t,\cdot)-\nu(t,\cdot)\Vert_{L^{2}(\vert\xi\vert\leq\delta_{1})} + \Vert w(t,\cdot)-\nu(t,\cdot)\Vert_{L^{2}(\delta_{1}\leq \vert\xi\vert\leq\delta)}\]
\[\leq C(\Vert u_{0}\Vert_{1} + \Vert u_{1}\Vert_{1,1})t^{-\frac{n}{2}} + \Vert w(t,\cdot)\Vert_{L^{2}(\delta_{1}\leq \vert\xi\vert\leq\delta)} + \Vert \nu(t,\cdot)\Vert_{L^{2}(\delta_{1}\leq \vert\xi\vert\leq\delta)}\]
\[\leq C(\Vert u_{0}\Vert_{1} + \Vert u_{1}\Vert_{1,1})t^{-\frac{n}{2}} + Ce^{-\alpha t}(\Vert u_{0}\Vert + \Vert u_{1}\Vert) + \vert P_{1}\vert e^{-\alpha t},\]
which implies the desired first estimate of Theorem 1.3.\\


Let us make sure the high frequency estimates for $w(t,\xi)$ directly. This seems easier because the characteristic roots are complex-valued on this zone.
\vspace{0.3cm}

\underline{{(3)}\,High frequency estimate for $w(t,\xi)$ when $\mu>2$}.
\\ 
On the zone of high  frequency the characteristic roots are complex-valued according to \eqref{roots-real-com}. Then the solution in the Fourier space given by \eqref{ike-19-1} can be rewritten for $|\xi|>\delta$ as
\begin{align}\label{roots-real-com-1}
w(t,\xi)& = (1+r)^{-\frac{t\mu}{2}}\frac{2}{r\sqrt{4-(\frac{\mu\log(1+r)}{r})^{2}}}\sin\big(\frac{r\sqrt{4-(\frac{\mu\log(1+r)}{r})^{2}}}{2}t\big)\;w_1(\xi) \nonumber\\
&+  (1+r)^{-\frac{t\mu}{2}} \cos\big(\frac{r\sqrt{4-(\frac{\mu\log(1+r)}{r})^{2}}}{2}t \big)\;w_0(\xi)  \\
& + (1+r)^{-\frac{t\mu}{2}}  \frac{ \mu \log(1+r)}{r\sqrt{4-(\frac{\mu\log(1+r)}{r}})^2}
\sin\big(\frac{r\sqrt{4-(\frac{\mu\log(1+r)}{r})^{2}}}{2}t \big) \;w_0(\xi) ,
\nonumber 
\end{align}
where we have used that $\exp(-\frac{t\mu}{2}\log(1+r))= (1+r)^{-\frac{t\mu}{2}}$.

Now, it is necessary to estimate the  three terms in the right hand side  of \eqref{roots-real-com-1}.

The first estimate is to 
\begin{align}\label{int-1}
\int_{r=|\xi|>\delta}&\Big| (1+r)^{-\frac{t\mu}{2}}\frac{2}{r\sqrt{4-(\frac{\mu\log(1+r)}{r})^{2}}}\sin\big(\frac{r\sqrt{4-(\frac{\mu\log(1+r)}{r})^{2}}}{2}t\big)\;w_1(\xi)\Big|^2d\xi\nonumber\\
&\leq t^{2}\int_{|\xi|>\delta} (1+r)^{-{t\mu}}\;|w_1(\xi)|^2d\xi\nonumber \\
&\leq \omega_n\|w_1\|_{\infty}^2 t^{2}\int_{\delta}^{\infty}(1+r)^{-t\mu}r^{n-1}dr \nonumber \\
& \leq C\|u_1\|_{1}^2 t^{2}e^{-\beta_1t} \leq C\|u_1\|_{1}^2 e^{-\beta t},  \quad t \gg 1,
\end{align}
where  $C$ and $\beta$ is positive constants depending  on $n$ and $\mu$ and $\beta_1$ is some positive constant given by  Lemma \ref{lemma2}.
We have just used the estimate $|\sin\theta| \leq \vert\theta\vert$ for all $\theta \in {\bf R}$.

Similarly, to get the second estimate,  using the fact $|\cos\theta| \leq 1$ for all $\theta$, one has
\begin{align}\label{int-2}
\int_{r=|\xi|>\delta}& \Big| (1+r)^{-\frac{t\mu}{2}} \cos\big(\frac{r\sqrt{4-(\frac{\mu\log(1+r)}{r})^{2}}}{2}t \big)\;w_0(\xi)\Big|^2 d\xi
\nonumber \\
& \leq \int_{|\xi|>\delta}  (1+r)^{-{t\mu}}\;|w_0(\xi)|^2 d\xi
\nonumber\\
& \leq C\|u_0\|_{1}^2 e^{-\beta_1t},  \quad t \gg 1.
\end{align}

To estimate the third integral in the right hand side of \eqref{roots-real-com-1} we  proceed similarly as in the estimate to the first integral.
Then one has
\begin{align}\label{int-3}
 \int_{|\xi| >\delta}&\Big|(1+r)^{-\frac{t\mu}{2}}  \frac{ \mu \log(1+r)}{r\sqrt{4-(\frac{\mu\log(1+r)}{r}})^2}
\sin\big(\frac{r\sqrt{4-(\frac{\mu\log(1+r)}{r})^{2}}}{2}t \big) \;w_0(\xi)  \Big|^2d\xi\nonumber\\
 & \leq \frac{t^{2}}{4}\int_{|\xi| >\delta}(1+r)^{-t\mu}  \mu^2 \log^2(1+r)
 \;|w_0(\xi)|^2d\xi\nonumber\\
 & \leq \frac{t^{2}}{4}\int_{|\xi| >\delta}(1+r)^{-t\mu}  \mu^2 r^2
  \;|w_0(\xi)|^2d\xi\nonumber
  \\
 & \leq \frac{t^{2}}{4}\mu^2 \omega_n \|u_0\|_{1}^2\int_{\delta}^{\infty}(1+r)^{-t\mu}   r^{n+1} dr  \nonumber\\
 & \leq Ct^{2}\|u_0\|_{1}^2 e^{-\beta_2t} \leq C\|u_0\|_{1}^2 e^{-\beta t},  \quad t \gg 1,
\end{align}
for some $\beta_2$ and $\beta$ positive constants.

By combining estimates \eqref{int-1}, \eqref{int-2} and \eqref{int-3} we arrive at  the exponential estimate in the high frequency zone. Note that one can replace the integral region $\vert\xi\vert < \delta$ with $\vert\xi\vert \leq \delta$ in the result below since all derived final estimates are independent from $\delta > 0$. In fact, $\delta > 0$ is a fixed constant that depends only on $\mu > 0$ (see \eqref{R-00}).  
\begin{lem}\label{lem4.2} Let $\mu > 2$, and let $u_{0} \in L^{1}({\bf R}^{n})$ and $u_{1} \in L^{1}({\bf R}^{n})$. Then, it holds that
\begin{align}\label{high-mu-large-1}
\int_{|\xi| \geq \delta} |w(t,\xi)|^2d\xi \leq C_{n,\mu} \big(\|u_0\|_{1}^2 + \|u_1\|_{1}^2 \big)\; e^{-\beta t}, \quad t\gg 1,
\end{align}
where $C=C_{n,\mu}$ is a positive constant and $\beta >0$ is a constant.
\end{lem}

Finally, the second inequality of Theorem \ref{theorem-3} is a direct consequence of Lemma 3.3.

\section{Optimal $L^{2}$-estimates of solutions: the case $ \mu >2$}

In this section we shall try to get the optimal estimates of the solutions to problem (1.1)-(1.2) in terms of the $L^{2}$-norm. We treat the effective damping case $\mu > 2$.\\

In order to get such estimates of the quantity $\Vert u(t,\cdot)\Vert = C\Vert w(t,\xi)\Vert$ with some $C > 0$, we first note the following inequalities.

\[\Vert w(t,\cdot)\Vert^{2} \geq \int_{\vert\xi\vert \leq \delta_{1}}\vert w(t,\xi)-\nu(t,\xi) + \nu(t,\xi)\vert^{2}d\xi\]
\begin{equation}\label{CR-1}
\geq \frac{1}{2}\int_{\vert\xi\vert \leq \delta_{1}}\vert \nu(t,\xi)\vert^{2}d\xi - \int_{\vert\xi\vert \leq \delta_{1}}\vert w(t,\xi) - \nu(t,\xi)\vert^{2}d\xi,
\end{equation} 
because of $\vert a+b\vert^{2} \geq \frac{1}{2}\vert a\vert^{2}-\vert b\vert^{2}$ ($a,b \in {\bf C}$), and 
\[\Vert w(t,\cdot)\Vert^{2} = \left(\int_{\vert\xi\vert \leq \delta_{1}} +  \int_{\delta_{1} \leq \vert\xi\vert \leq \delta} +  \int_{\delta \leq \vert\xi\vert}\right) \vert w(t,\xi)\vert^{2}d\xi \]
\[\leq 2\int_{\vert\xi\vert \leq \delta_{1}}\vert w(t,\xi) - \nu(t,\xi)\vert^{2}d\xi + 2\int_{\vert\xi\vert \leq \delta_{1}}\vert \nu(t,\xi)\vert^{2}d\xi \]
\begin{equation}\label{CR-2}
+ \int_{\delta_{1} \leq \vert\xi\vert \leq \delta}\vert w(t,\xi)\vert^{2}d\xi + \int_{\delta \leq \vert\xi\vert}\vert w(t,\xi)\vert^{2}d\xi,
\end{equation} 
where $\delta > 0$ and $\delta_{1} > 0$ are constants defined in \eqref{roots-real-com} and Lemma \ref{ike-705}, respectively. So, because of Lemma \ref{lem4.1}, \eqref{ike-955} and \eqref{high-mu-large-1}, in order to prove Theorem 1.4 it suffices to consider only one factors
\[I_{0}(t) := \int_{\vert\xi\vert \leq \delta_{1}}\vert \nu(t,\xi)\vert^{2}d\xi.\]
We will get upper and lower bound estimates for $I_{0}(t)$ as $t \to \infty$.\\

We first get upper bound estimate.

\underline{Upper bound estimate. }

\[I_{0}(t) = P_{1}^{2}\int_{\vert\xi\vert \leq \delta_{1}}\frac{\vert e^{\lambda_{+}t} - e^{\lambda_{-}t}  \vert^{2}}{\vert\lambda_{+}-\lambda_{-}  \vert^{2}}d\xi = P_{1}^{2}t^{2}\int_{\vert\xi\vert \leq \delta_{1}}e^{2\lambda_{+}t} \frac{\vert 1 - e^{-(\lambda_{+}-\lambda_{-})t}\vert^{2}}{\vert(\lambda_{+}-\lambda_{-})t\vert^{2}}d\xi.\]

Then,  from \eqref{CR-7} one has
\[I_{0}(t) \leq P_{1}^{2}t^{2}\kappa_{0}^{2}\int_{\vert\xi\vert \leq \delta_{1}}e^{2\lambda_{+}t}d\xi = P_{1}^{2}t^{2}\kappa_{0}^{2}\omega_{n}\int_{0}^{\delta_{1}} e^{2\lambda_{+}t}r^{n-1}dr\]
\[\leq P_{1}^{2}t^{2}\kappa_{0}^{2}\omega_{n}\int_{0}^{\delta_{1}} e^{-2drt}r^{n-1}dr \leq CP_{1}^{2}t^{2}\kappa_{0}^{2}\omega_{n}t^{-n},\]
which implies the desired upper bound estimate
\begin{equation}\label{CR-4}
I_{0}(t) \leq CP_{1}^{2}t^{2-n},\quad t \gg 1, 
\end{equation} 
where one has just used (2) of Lemma \ref{ike-705}.

\underline{Lower bound estimate.}

Next, we get the lower bound estimate for $I_{0}(t)$. Let $\delta_{1} > 0$ and $\delta_{2} > 0$ be constants defined in Lemma \ref{ike-705} and \eqref{CR-3}, respectively, and choose $t > 0$ sufficiently large to satisfy
\[e^{\frac{\delta_{2}}{dt}} -1\leq \delta_{1}.\]
If we consider $r \in [0,e^{\frac{\delta_{2}}{dt}} -1]$, then $r \in [0,\delta_{1}]$ and 
\[d\log(1+r)t \leq \delta_{2},\]
where $d > 0$ is a constant defined in Lemma \ref{ike-705}. This implies, because (1) of Lemma \ref{ike-705}, 
\[(\lambda_{+}-\lambda_{-})t \leq d\log(1+r) t \leq \delta_{2}, \quad t \gg 1.\]
Thus, one can apply \eqref{CR-3} of Section 3 to get
\begin{equation}\label{ike-950}
\gamma_{0} \leq \frac{1-e^{-(\lambda_{+}-\lambda_{-})t}}{(\lambda_{+}-\lambda_{-})t}, \quad r \in (0,e^{\frac{\delta_{2}}{dt}} -1], \quad t\gg1.
\end{equation}
Based on \eqref{ike-950} and (2) of Lemma \ref{ike-705} one can estimate $I_{0}(t)$ as follows:
\[I_{0}(t) = P_{1}^{2}t^{2}\int_{\vert\xi\vert\leq\delta_{1}}e^{2t\lambda_{+}}\left\vert \frac{1-e^{-(\lambda_{+}-\lambda_{-})t}}{(\lambda_{+}-\lambda_{-})t} \right\vert^{2}d\xi\]
\[\geq  P_{1}^{2}t^{2}\gamma_{0}^{2}\int_{\vert\xi\vert\leq e^{\frac{\delta_{2}}{dt}} -1}e^{2t\lambda_{+}}d\xi \geq  P_{1}^{2}t^{2}\gamma_{0}^{2}\int_{\vert\xi\vert\leq  e^{\frac{\delta_{2}}{dt}} -1}e^{-2ctr}d\xi\]
\[\geq  P_{1}^{2}t^{2}\gamma_{0}^{2}\omega_{n}\int_{0}^{e^{\frac{\delta_{2}}{dt}} -1}e^{-2ctr}r^{n-1}dr, \quad t\gg 1,\]
where $c$ is a positive constant defined in (2) of Lemma \ref{ike-705}.

Since $e^{x} - 1 > x$ for $x \geq 0$ one can get the estimate
\[I_{0}(t) \geq  P_{1}^{2}t^{2}\gamma_{0}^{2}\omega_{n}\int_{0}^{\frac{\delta_{2}}{dt}}e^{-2ctr}r^{n-1}dr \geq P_{1}^{2}t^{2}\gamma_{0}^{2}\omega_{n}e^{-2ct(\frac{\delta_{2}}{dt})}\int_{0}^{\frac{\delta_{2}}{dt}}r^{n-1}dr\]
\[=  P_{1}^{2}t^{2}\gamma_{0}^{2}\omega_{n}\frac{e^{-2c\delta_{2}/d}}{n}\frac{\delta_{2}^{n}}{d^{n}}t^{-n},\]
which implies the desired lower bound estimate: 
\begin{equation}\label{CR-8}
I_{0}(t) \geq CP_{1}^{2}t^{2-n},\quad t \gg 1. 
\end{equation} 

The proof of Theorem 1.4 is a direct consequence of \eqref{CR-1}, \eqref{CR-2}, \eqref{CR-4}, \eqref{CR-8} as is already mentioned.

\begin{rem}{\rm From the proof above one may realize that one can not treat two factors $P_{1}\frac{e^{\lambda_{+}t}}{\sqrt{\mu^{2}\log^{2}(1+r)-4r^{2}}}$ and $P_{1}\frac{e^{\lambda_{-}t}}{\sqrt{\mu^{2}\log^{2}(1+r)-4r^{2}}}$ of the leading term $\nu(t,\xi)$ separately. In this sense, one can observe an effective role of the "double diffusion" structure of the solution itself as is pointed out in \cite{DE} and \cite{PCI}.}
\end{rem}


\section{Optimal $L^{2}$-estimates of solutions:  case $0< \mu < 2$}

In this section we shall prove Theorem 1.2. To do that we treat the non-effective damping case $\mu \in (0,2)$ where the characteristic roots associated with the problem
\eqref{F-eqn}--\eqref{F-initial} are complex-valued for all $\xi \in {\bf R}^{n}$.

For the estimates needed in this section, one uses the following facts that
\begin{align}\label{ike-100}
L := \sup_{\theta \ne 0}\left\vert\frac{\sin\theta}{\theta}\right\vert < +\infty,
\end{align}
and there exists a real number $\delta_{0} \in (0,1)$ such that for all $\theta \in (0,\delta_{0}]$
\begin{equation}\label{ike-101}
\left\vert\frac{\sin\theta}{\theta}\right\vert \geq \frac{1}{2}.
\end{equation}
As in the same concept discussed in Section 4, for our purpose it suffices to get upper and lower bound estimates for the quantity 
\[K_{n}(t) := \int_{\vert\xi\vert \leq 1}\vert\chi(t,\xi)\vert^{2}d\xi,\] 
where
\begin{eqnarray}
\chi(t,\xi)  := (1+r)^{-\frac{t\mu}{2}}\frac{2}{r\sqrt{4-(\frac{\mu\log(1+r)}{r})^{2}}}\sin(\gamma tr)P_{1},
\end{eqnarray}
and $\gamma = \displaystyle{\frac{\sqrt{4-\mu^{2}}}{2}} > 0$.
\noindent
Note that a contribution from the high frequency estimates in $\{\vert\xi\vert \geq 1\}$ is exponentially small, and in fact, the following exponential decay estimate is true:
\begin{align}\label{high-freq-nu}
\int_{|\xi| \geq 1}\vert \chi(t,\xi)\vert^2 d\xi \leq P_1^2 t^{2}L^{2}\int_{|\xi| \geq 1} (1+r)^{-\mu t}d\xi \leq C_{n}P_1^2 t^{2}e^{-\alpha t}, \quad t > 0,
\end{align}
because of Lemma 2.2, where the constants $C_{n} > 0$ and $\alpha > 0$ depend on $n$, and $\mu \in (0,2)$. 

We first get the lower bound estimate for $K_{n}(t)$. For this purpose, we prepare one facts that since
\[\lim_{r \to +0}\frac{1}{4-(\frac{\mu\log(1+r)}{r})^{2}} = \frac{1}{4\gamma^{2}},\]
there exists a constant $\rho_{0} \in (0,1)$ such that for all $r \leq \rho_{0}$, one has
\begin{equation}\label{ike-120}
\frac{1}{8\gamma^{2}} \leq \frac{1}{4-(\frac{\mu\log(1+r)}{r})^{2}} \leq \frac{1}{2\gamma^{2}}.
\end{equation}
Now, take $t > 0$ sufficiently large such that $\displaystyle{\frac{\delta_{0}}{\gamma t}} < \rho_{0}$. This implies $\gamma tr \leq \delta_{0}$ and $r < \rho_{0}$ if $r \leq \displaystyle{\frac{\delta_{0}}{\gamma t}}$, which brings to the useful inequalities \eqref{ike-120} and 
\begin{equation}\label{ike-107}
\left\vert\frac{\sin(\gamma tr)}{\gamma tr }\right\vert \geq \frac{1}{2},
\end{equation}
where one has just used \eqref{ike-101}. Then, from \eqref{ike-120} and \eqref{ike-107} one has a series of inequalities:
\[K_{n}(t) \geq 4P_{1}^{2}\int_{\vert\xi\vert \leq \rho_{0}}(1+r)^{-\mu t}\frac{\sin^{2}(\gamma rt)}{r^{2}\left(4-(\frac{\mu\log(1+r)}{r})^{2}\right)}d\xi \]
\[= 4t^{2}P_{1}^{2}\gamma^{2}\int_{\vert\xi\vert \leq \rho_{0}}(1+r)^{-\mu t}\frac{\sin^{2}(\gamma rt)}{t^{2}\gamma^{2}r^{2}}\left(\frac{1}{(4-(\frac{\mu\log(1+r)}{r})^{2}}\right)d\xi \]
\[\geq 4t^{2}P_{1}^{2}\gamma^{2}\int_{\vert\xi\vert \leq \frac{\delta_{0}}{\gamma t}}(1+r)^{-\mu t}\frac{\sin^{2}(\gamma rt)}{t^{2}\gamma^{2}r^{2}}\left(\frac{1}{(4-(\frac{\mu\log(1+r)}{r})^{2}}\right)d\xi \]
\[\geq t^{2}\frac{P_{1}^{2}}{8}\int_{\vert\xi\vert \leq  \frac{\delta_{0}}{\gamma t}}(1+r)^{-\mu t} d\xi = t^{2}\frac{P_{1}^{2}}{8}\omega_{n}\int_{0}^{ \frac{\delta_{0}}{\gamma t}}(1+r)^{-\mu t}r^{n-1}dr\]
\[\geq t^{2}\frac{P_{1}^{2}}{8}\omega_{n}(1+\frac{\delta_{0}}{\gamma t})^{-\mu t}\frac{1}{n}(\frac{\delta_{0}}{\gamma})^{n}t^{-n}\]
\begin{equation}\label{ike-102}
\geq \frac{\omega_{n}P_{1}^{2}}{16n}(\frac{\delta_{0}}{\gamma})^{n}e^{-\frac{\mu\delta_{0}}{\gamma}}t^{2-n} =: C_{n}t^{2-n},\quad t \gg 1,
\end{equation}
where one has just relied on the fact that
\[\lim_{t \to \infty}(1+\frac{\delta_{0}}{\gamma t})^{-\mu t} = e^{-\frac{\mu\delta_{0}}{\gamma}}.\]

While,  one can get the upper bound estimates for $K_{n}(t)$. Indeed, one has
\[K_{n}(t) = \left(\int_{\vert\xi\vert\leq \rho_{0}} +\int_{1 \geq \vert\xi\vert \geq \rho_{0}}\right)\vert\chi(t,\xi)\vert^{2}d\xi =: K_{n,1}(t) + K_{n,2}(t).\]
$K_{n,1}(t)$ can be estimated as follows from \eqref{ike-100}, \eqref{ike-120} and Lemma 2.1:
\[K_{n,1}(t) \leq 4\gamma^{2}P_{1}^{2}t^{2}L^{2}\int_{\vert\xi\vert \leq \rho_{0}}(1+r)^{-\mu t}\frac{1}{4-(\frac{\mu\log(1+r)}{r})^{2}}d\xi\]
\begin{equation}\label{ike-103}
\leq 2P_{1}^{2}t^{2}L^{2}\int_{\vert\xi\vert \leq \rho_{0}}(1+r)^{-\mu t}d\xi \leq D_{n}P_{1}^{2}t^{2-n}
\end{equation}
for large $t > 1$ and some constant $D_{n} > 0$. Furthermore, for the estimate of $K_{n,2}(t)$ one can proceed the same computation as in \eqref{high-freq-nu} to get
\begin{align}\label{ike-111}
K_{n,2}(t) \leq \int_{|\xi| \geq \rho_{0}}\vert \chi(t,\xi)\vert^2 d\xi \leq C_{n}P_1^2 t^{2}e^{-\alpha t}, \quad t > 0.
\end{align}

By using \eqref{ike-102}, \eqref{ike-111} and \eqref{ike-103} one can obtain the following lemma:
\begin{lem}\label{ike-37}Let $n \geq 1$ and $\mu \in (0,2)$. Than it holds that
\[C_{n}^{2}P_{1}^{2}t^{2-n} \leq \int_{\vert\xi\vert \leq 1}\vert\chi(t,\xi)\vert^{2}d\xi \leq D_{n}^{2}P_{1}^{2}t^{2-n}\quad (t \gg 1),\]
where $C_{n} > 0$ and $D_{n} > 0$ are generous constants depending on $n$ and $\mu \in (0,2)$. 
\end{lem}

To end this section, based on the inequalities

\[\int_{|\xi| \leq 1}\vert \chi(t,\xi) \vert^2d\xi - \int_{{\bf R}^n}\vert w(t,\xi)-\chi(t,\xi)\vert^2d\xi \leq \int_{{\bf R}^n}\vert w(t,\xi)\vert^2d\xi \]
\begin{equation}\label{main-ineq}
\leq \int_{{\bf R}^n}\vert w(t,\xi)-\chi(t,\xi)\vert^2d\xi + \int_{ |\xi| \leq 1}\vert \chi(t,\xi) \vert^2d\xi + \int_{ |\xi| \geq 1}\vert \chi(t,\xi) \vert^2d\xi, 
\end{equation}
because of Theorem 1.1, estimates \eqref{high-freq-nu} and Lemma \ref{ike-37} combined with the Plancherel Theorem, one can prove Theorem \ref{theorem-2-1}.\\


\section{Optimal $L^{2}$-estimates of solutions:  case $\mu = 2$}

In this section, we prove Theorem \ref{theorem-ike-2} based on Theorem 1.5.
For this purpose, one shall consider the critical case of $\mu = 2$. In this case, from Theorem 1.5 one can see that the leading term is 
\[\nu(t,\xi) := (1+t)^{-t}\frac{\sin\left(t\sqrt{r^{2}-\log^{2}(1+r)}\right)}{\sqrt{r^{2}-\log^{2}(1+r)}}P_{1}\quad r > 0.\]  

It suffices to estimate the quantity: for $P_{1} \ne 0$;
\begin{equation}\label{ike-500}
\frac{1}{P_{1}^{2}}\int_{\vert\xi\vert \leq \eta}\vert\nu(t,\xi)\vert^{2}d\xi
\end{equation}
for small $\eta > 0$.

Now, note that since
\[\lim_{r \to +0}\frac{r^{2}-\log^{2}(1+r)}{r^{3}} = 1,\]
\[\lim_{r \to +0}\frac{r-\frac{\log(1+r)}{1+r}}{r^{2}} = \frac{3}{2},\]
one can assume that
\begin{equation}\label{ike-200}
r^{2}-\log^{2}(1+r) \sim r^{3} \quad r \to +0,
\end{equation}
\begin{equation}\label{ike-201}
r-\frac{\log(1+r)}{1+r} \sim r^{2} \quad r \to +0.
\end{equation}
One uses \eqref{ike-200}, \eqref{ike-201} and a polar coordinate transform to estimate the essential part of \eqref{ike-500}:
\[I_{n}(t) := \int_{0}^{\eta}(1+r)^{-2t}r^{n-1}\frac{\sin^{2}(t\sqrt{f(r)})}{(\sqrt{f(r)})^{2}}dr ,\]
where $f(r) := r^{2}-\log^{2}(1+r)$ for small $\eta > 0$. One applies the change of variable as $w = \sqrt{f(r)}$. Then, it is true that
\[dw = \frac{r-\frac{\log(1+r)}{1+r}}{\sqrt{f(r)}}dr.\]
It should be noticed that  $w \sim \sqrt{r^{3}} = r^{\frac{3}{2}}$, so that $r \sim w^{\frac{2}{3}}$ when $r \to +0$. Then, one has a series of inequalities: for small $\eta > 0$ it holds that
\[I_{n}(t) \sim \int_{0}^{\eta}(1+w^{\frac{2}{3}})^{-2t}w^{\frac{2n-6}{3}}\frac{\sin^{2}(tw)}{w}dw = \int_{0}^{\eta}(1+w^{\frac{2}{3}})^{-2t}w^{\frac{2n-3}{3}}\frac{\sin^{2}(tw)}{w^{2}}dw =:\tilde{I}_{n}(t).\]

\underline{Lower bound estimate for $I_{n}(t) \sim \tilde{I}_{n}(t)$}:\\
By a change of variable by $\sigma = w^{\frac{2}{3}}$ one has
\[\tilde{I}_{n}(t) = \frac{3}{2}\int_{0}^{\eta^{\frac{2}{3}}}(1+\sigma)^{-2t}\sigma^{n-4}\sin^{2}(t\sigma^{\frac{3}{2}})d\sigma.\]
Since $\log(1+\sigma) \sim \sigma$ as $\sigma \to 0$, one finds that $(1+\sigma)^{-2t} \sim e^{-2t\sigma}$ for small $\sigma > 0$. Thus, one can asume
\begin{equation}\label{ike-600}
\tilde{I}_{n}(t) \sim \int_{0}^{\eta^{\frac{2}{3}}}e^{-2t\sigma}\sigma^{n-4}\sin^{2}(t\sigma^{\frac{3}{2}})d\sigma
\end{equation}
\[= t^{2}\int_{0}^{\delta_{0}}e^{-2t\sigma}\sigma^{n-1}\left\vert\frac{\sin(t\sigma^{\frac{3}{2}})}{t\sigma^{\frac{3}{2}}}\right\vert^{2}d\sigma\] 
for small $\delta_{0} > 0$ defined in \eqref{ike-101}. Thus, by a change of variable $t\sigma = r$ one has
\[\tilde{I}_{n}(t) \geq C\frac{t^{2}}{4}\int_{0}^{\delta_{0}^{\frac{2}{3}}t^{-\frac{2}{3}}}e^{-2t\sigma}\sigma^{n-1}d\sigma\]
\[= C\frac{t^{2-n}}{4}\int_{0}^{At^{\frac{1}{3}}}e^{-2r}r^{n-1}dr \geq C\frac{t^{2-n}}{4}\int_{0}^{At^{\frac{1}{3}}}e^{-2r}r^{n-1}dr\]
\[\geq C\frac{t^{2-n}}{4}\int_{0}^{1}e^{-2r}r^{n-1}dr,\]
for large $t > 1$ satisfying $At^{\frac{1}{3}} > 1$, where $A := \delta_{0}^{\frac{2}{3}}$ and $C > 0$ is a constant brought by equivalence \eqref{ike-600}. This implies that there exists a generous constant $C_{n} > 0$ such that
\begin{equation}\label{ike-202}
I_{n}(t) \sim \tilde{I}_{n}(t) \geq C_{n}t^{2-n}, \quad t \gg 1.
\end{equation}

\underline{Upper bound estimate $I_{n}(t) \sim \tilde{I}_{n}(t)$}:\\
Because of \eqref{ike-100} one can estimate
\[\tilde{I}_{n}(t) \leq L^{2}t^{2}\int_{0}^{\eta}(1+w^{\frac{2}{3}})^{-2t}w^{\frac{2n-3}{3}}dw\]
\[= L^{2}t^{2}\int_{0}^{\eta^{2/3}}(1+\sigma)^{-2t}\sigma^{n-1}dw\]
\begin{equation}\label{ike-203}
\leq CL^{2}t^{2}t^{-(n-1+1)} = CL^{2}t^{2-n}
\end{equation}
with some constant $C > 0$ because of Lemma 2.1 with $p := n-1$.

A rest part of proof of Theorem \ref{theorem-ike-2} can be done as in Section 5 based on the inequality \eqref{main-ineq} with $\chi(t,\xi)$ replaced by $\nu(t,\xi)$ above.

\par
\vspace{0.5cm}
\noindent{\em Acknowledgement.}
\smallskip
The work of the first author (R. C. Char\~ao) was partially supported by Print/Capes - Process 88881.310536/2018-00 and the work of 
 the second author (R. Ikehata) was supported in part by Grant-in-Aid for Scientific Research (C)20K03682  of JSPS.


\end{document}